\magnification=1200
\outer\def\beginsection#1\par{\vskip0pt plus.1\vsize\penalty-250\vskip0pt plus -.1\vsize\bigskip\vskip\parskip\message{#1}\leftline{\bf#1}\nobreak\smallskip\noindent}
\def\cov{\hat V(Q,\,R)}
\def\vec{V(Q,\,R)}
\def\ux{{\underline x}}
\def\uz{{\underline z}}
\def\upar{{\underline \partial}}
\def\xe{x_1^{e_1}\cdots\,x_n^{e_n}}

\def\de{\partial_n^{e_n}\cdots\,\partial_1^{e_1}}
\def\df{\partial_n^{f_n}\cdots\,\partial_1^{f_1}}

\def\sqr#1#2{{\vcenter{\vbox{\hrule height.#2pt\hbox{\vrule width.#2pt height#1pt \kern#1pt \vrule width.#2pt}\hrule height.#2pt}}}}
\def\square{\mathchoice\sqr64\sqr64\sqr{2.1}3\sqr{1.5}3}

\def\ut{{\underline t}}

\def\te{t_1^{e_1}\cdots\,t_n^{e_n}}
\def\cuno{_{\cdot1}}
\def\unoc{_{1\cdot}\!}
\def\cove{\cov^{\rm ext}}
\def\vece{\vec^{\rm ext}}

\def\isigmagamma{I''_{(\sigma,\,\gamma)}}

\def\id{{\rm id}}

\centerline{\bf ``On the braided Fourier transform in the $n$-dimensional quantum space''}\vskip0.3cm
{\sl 
\centerline{ Giovanna Carnovale\footnote*{\rm This research is financially supported by the NWO-Stichting SWON, under project number 610-06-100}}

\centerline{ Mathematisch Instituut}

\centerline{P.O.Box 80.010,  \ 3508 TA }

\centerline{Utrecht, The Netherlands}

\centerline{email address: carnoval@math.uu.nl}}\vskip0.3cm 

{\noindent{\bf Abstract:} We work out in detail a theory of integrability on the braided covector Hopf algebra and the braided vector Hopf algebra of type $A_n$ introduced in [Ma] and [KeMa].  Starting with a definition of braided Fourier transform very similar to that in [KeMa] we obtain $n$-dimensional analogous results to those in [Koo] expressing the correspondence between products of the $q^2$-Gaussian $g_{q^2}(\ux)$ times monomials, and products of the $q^2$-Gaussian $G_{q^2}(\upar)$ times $q^2$-Hermite polynomials under the transform. We invert the correspondence by finding a suitable inversion, different from that in [KeMa]. We show that with this transforms, whenever $n\ge2$, the Plancherel measure will depend on the parity of the power series that we are transforming.} \smallbreak
 
\beginsection1.0. Introduction.

Recently Majid ([Ma], and references given there) has defined a generalisation of the concept of
Hopf algebra, namely ``braided Hopf algebras''. Hopf superalgebras and genuine Hopf algebras are examples of this objects, but there are more examples,
associated to quantum groups, since the category of representations of a quasitriangular Hopf algebra is braided. This is also the subject of [HH] with different terminology.

\noindent Kempf and Majid introduced in [KeMa] an integration theory for
a class of braided Hopf algebras arising from  matrix solutions of the quantum Yang-Baxter equation, as
``braided covector algebras''(see also [Ma]). They used this theory to define a formal braided Fourier transform  and its inverse on these algebras. In their paper they also present the case of the braided line as an example, and the $n$-dimensional case in less details. 

\noindent The main problem in their theory is that it is very difficult to
find an explicit integral that behaves well enough. They provide powerful general results in a theoretical way, but the description in specific cases is often hard to
 handle. Besides, they do not provide a definition
of convergence of an integral nor do they say in their article when an object is called integrable.

\noindent The purpose of this paper is to work out as far as we can the example
of the $n$-dimensional quantum space of type $A_n$ viewed as a braided Hopf algebra. We will provide a definition of an integral (similar to that in [KeMa]) accompanied by results on convergence of the integral based on extensions of representations of the braided covector and vector algebras. We define three different concepts of integrability, and provide examples and counterexamples. Using this facts, we show, similarly to [KeMa] that the given integral is invariant under translation, and we use it for the definition of two types of Fourier transform. Our transforms will be slightly different from that of [KeMa]. One of them looks more like that
of [Ch] since the integral does not have trivial braiding with elements in the braided Hopf algebra, and because the braided antipode appears in the definition. We also took inspiration from [Koo], who studied  also an analogue of the Fourier transform for the case of the braided line, although his transform goes from an algebra to itself, while we are looking for a transform going from an algebra to its braided dual, as in [KeMa], [OR] and [Ch].

\noindent We find an $n$-dimensional analogue of the correspondence between products of $q^2$-Gaussians times monomials and other $q^2$-Gaussians times $q^2$-Hermite polynomials, similar to the results for the braided line in [Koo]. 
\noindent We give also  inverses for our transforms that invert the correspondence mentioned above, similarly to what appears in [Koo]. The main tool for this inversion formula is the symmetry between the braided vector algebra and the braided covector algebra. Kempf and Majid had already defined an inversion formula in their article, but they used properties that our integral does not have. Other inversion formulae for the braided line are to be found in [OR] where the case of distributions is also treated.

\noindent One of the most interesting results is that whenever $n\ge2$, there is a loss of symmetry, so that the Plancherel measure will no longer be the same in the whole space. Indeed there is a action of ${\bf Z_2}^n$ associted to the parity of the power series we are working with, and the Plancherel measure will be constant only on the subspaces of power series with constant parity. Therefore the transforms we define can also be seen as sine and cosine transform. A phenomenon of break of simmetry for $q$-integrals was also noted  in [CHMW], where the authors were defining a calculus associated to a $q$-deformed Heisenberg Algebra.

\noindent Other definitions of analogues of the Fourier, on genuine Hopf algebras, quantum spaces, or commutative algebras appeared before [KeMa] in [MMNNU], [KooSw] and [LyuMa].

\beginsection1.1. Notation

\noindent In this paper a complex algebra has, unless otherwise stated, always a unit. $q$ is a real number $0<q<1$. \noindent For a positive integer $m$, and for any $q\not=1$ we write $[m]_q={{q^{m}-1}\over{q-1}}$ and $[m]_{q}!=\prod_{j=1}^m[j]_q$.

\noindent We will also use the $q-$shifted factorial. For any $a\in{\bf R}$ and for any $k\in{\bf Z}_{\ge0}$, we will then put $(a;q)_k=\prod_{l=0}^{k-1}(1-aq^l)$. If $0<q<1$ we will also write $(a;q)_{\infty}=\lim_{k\to\infty}(a;q)_k$ and for $r$ real numbers $a_1,\ldots,\,a_r$ we will put $(a_1,\ldots,a_r;q)_\infty=\prod_{j=1}^r(a_j;q)_\infty$. Finally, for $a\ge b$  with $a$ and $b$ both in ${\bf Z}_{\ge0}$ in we will use the $q-$binomial coefficient $\left[{a\atop b}\right]_q={{[a]_q!}\over{[b]_q![a-b]_q!}}={{(q;q)_a}\over{(q;q)_b(q;q)_{a-b}}}$ 

\noindent Whenever for any capital character $E$ we have a multi-index $E=(e_1,\,\ldots,\,e_n)$ we will put: $E_i=\sum_{j=1}^{i-1}e_j$ and $E^i=\sum_{j=i+1}^ne_j$. Hence $|E|=e_i+E_i+E^i$ for every $i$.

\noindent We identify the set $\{+,\,-\}$ with the field $\bf Z_2$, letting $+$ correspond to  $\bar 0$ and $-$ correspond to $\bar1$, so that 
$n$-tuples in $\{+,\,-\}^n$ can be identified with vectors in ${\bf Z_2}^n$. By means of this identification we define the map $A\colon{\bf Z}^n\to{\bf Z}_2^n\to\{+,\,-\}^n$ by reducing modulo $2$ first. I.e. $A(b)=+$ if $b$ is even and $-$ otherwise. We will also denote by $B\colon\{+,\,-\}^n\to \{0,\,1\}^n\subset{\bf Z}^n$  the map sending each ``even'' entry to $0$ and each ``odd'' entry to $1$. 

\noindent Given an operator on the two-fold tensor product of an $n-$dimensional  vector space $V$, we identify this operator with the $n^2\times n^2$ matrix $R$ and we denote its entries by $R^{ab}_{cd}$ where $a,\, b$ are the row entries and $c,\,d$ are the column entries. For such an $R$,  the operator acting on the $p-$fold tensor product of $V$ (for $p\ge2$) as $R$ on the $i^{th}$ and $j^{th}$ component and as the identity elsewhere will be denoted by $R_{ij}$. For summation we will use Einstein convention.

\noindent For any vector space $V$, we will denote by $P$ the operator $P\colon V\otimes V\to V\otimes V$ sending $v\otimes w$ to $w\otimes v$.  

\beginsection2. Braided covector algebra and braided vector algebra.

\noindent We begin with the definition of braided bialgebras and Hopf algebras, and we give then the definition of the braided covector algebra
and braided vector algebra as given in [Ma].  

\proclaim Definition. A braided Hopf algebra over a field $K$ is an associative  algebra $A$ with multiplication $m$ and a coassociative coalgebra with comupltiplication $\Delta$ and counit $\varepsilon$ together with an invertible  linear map $\Phi\colon A\otimes A\to A\otimes A$ called braiding,  and a linear map $S\colon A\to A$ called braided antipode such that the following properties hold:   $$\eqalign{&\Phi(m\otimes\id)=(\id\otimes m)(\Phi\otimes \id)(\id\otimes\Phi)\cr
&\Phi(\id\otimes m)=(m\otimes \id)(\id\otimes \Phi)(\Phi\otimes\id)\cr
&(\id\otimes\Delta)\circ \Phi=(\Phi\otimes\id)(\id\otimes \Phi)(\Delta\otimes \id)\cr
&(\Delta\otimes \id)\circ \Phi=(\id\otimes \Phi)(\Phi\otimes \id)(\id\otimes\Delta)\cr
&1\varepsilon=m(\id\otimes S)\Delta=m(S\otimes \id)\Delta\cr
&\Delta\,m=(m\otimes m)(\id \otimes\Phi\otimes \id)(\Delta\otimes\Delta)\cr
&\varepsilon\circ m=\varepsilon\otimes \varepsilon\cr
&\Delta(1)=1\otimes 1\cr}$$ i.e. $\Delta$ is an algebra map $A\to A\otimes A$ where $A\otimes A$ has the product  $(m\otimes m)(\id \otimes\Phi\otimes \id)$ that is associative and coassociative because of the properties of the braiding, and $\varepsilon$ is an algebra map $A\to K$. If there is no braided antipode, $A$ is called ``braided bialgebra''.\par 

\noindent S. Majid who first introduced the concept, defines braided Hopf algebras as Hopf algebras living in a braided category. The definition we give here is equivalent to Majid's definition if the objects in the category are vector spaces. For the categorical approach and an extensive description of the properties and importance of such objects, see [Ma] and references given there. For a definition similar to the one appearing here, see [Du]. In the latter one can find a more general definition than here because the counit is not required to be an algebra homomorphism. [Du] shows that the two definitions coincide if and only if $\varepsilon$ is an algebra homomorphism and he describes a few properties and equivalences. In particular we recall that in a braided Hopf algebra there holds 
$$\eqalign{&\Phi(1\otimes a)=a\otimes 1,\quad\Phi(a\otimes 1)=1\otimes a\quad\forall a\,\in\,A\cr
&(\varepsilon\otimes \id)\Phi=\id\otimes\varepsilon,\quad(\id\otimes\varepsilon)\Phi=(\varepsilon\otimes \id)\cr
&\Phi=(m\otimes m)(S\otimes \Delta\otimes S)(\id\otimes m\otimes \id)(\Delta\otimes \Delta)\cr
&(\Phi\otimes \id)(\id\otimes \Phi)(\Phi\otimes\id)=(\id\otimes \Phi)(\Phi\otimes \id)(\id\otimes\Phi)\cr
&S(1)=1,\quad \varepsilon\circ S=\varepsilon\cr
&\Delta\,S=\Phi(S\otimes S)\Delta,\quad S\,m=m(S\otimes S)\Phi\cr
&\Phi(S\otimes \id)=(\id\otimes S)\Phi,\quad
\Phi(\id\otimes S)=(S\otimes \id)\Phi\cr}$$ In particular we observe that the operator $\Phi$ must satisfy the braid relation. On the other hand whenever we have an operator $T$ satisfying the braid relations, it is well known that the operator $R=TP$ satisfies the quantum Yang-Baxter equation in the form  $R_{12}R_{13}R_{23}=R_{23}R_{13}R_{12}$. Hence one can build the dual quasitriangular bialgebra $A(R)$ by means of the FRT construction (see for instance [FRT]). One can find in [Ma, Chapter 9] and references in there that the category of right comodules for $A(R)$ is braided, and that $R$
defines a braiding between every pair of right comodules. In particular for any two comodule algebras $U$ and $W$, there is an algebra structure on $U\otimes W$ where the product is defined by means of the braiding and the products in $U$ and $W$, such that $U\otimes W$ turns out to be again a comodule algebra for $A(R)$. This can be found in [Ma, Chapter 9]. Moreover, for a sufficiently well behaved braided category, on can always find a dual quasitriangular Hopf algebra $H$ such that the category we have is the category of right comodules  for $H$. 

\noindent We will work with the following example of braided Hopf algebras:
 
\proclaim Definition-Proposition [Ma, Theorem 10.2.1]. Let $R$ be an invertible $n^2\times n^2$ matrix obeying the quantum Yang-Baxter equation. Suppose also that there is another invertible $n^2\times n^2$ matrix $Q$ for which $Q_{12}R_{13}R_{23}=R_{23}R_{13}Q_{12}$
and $R_{12}R_{13}Q_{23}=Q_{23}R_{13}R_{12}$, and $(PR+1)(PQ-1)=0$. Then the ``braided covector algebra'' $\cov$ (resp. ``braided vector algebra'' $\vec$) defined by the generators $x_1,\,\ldots,\,x_n$
(resp. $\partial^1,\,\ldots\,\partial^n$) and the relations $x_ix_j=Q^{ab}_{ij}x_bx_a$ (resp. $\partial^i\partial^j=Q^{ij}_{ab}\partial^b\partial^a$) form a braided bialgebra with $\Delta(x_i)=x_i\otimes 1+1\otimes x_i$ (resp. $\Delta(\partial^i)=\partial^i\otimes 1+1\otimes \partial^i$), $\varepsilon(x_i)=0$ (resp. $\varepsilon(\partial^i)=0$) and braiding $\Phi(x_i\otimes x_j)=x_b\otimes x_aR^{ab}_{ij}$ (resp. $\Phi(\partial^i\otimes\partial^j)=R^{ij}_{ab}\partial^b\otimes\partial^a$). Moreover, if $PRPQ=PQPR$ we have a braided Hopf algebra with braided antipode
$S(x_i)=-x_i$ (resp. $S(\partial^i)=-\partial^i)$.\par

\noindent{\bf Proof:} We sketch the idea of the proof. The quantum Yang-Baxter equation for $R$ ensures that the braid relations hold for $\Phi$ on the generators. The relations involving $Q_{ij}$ and $R_{kl}$ ensure that $\Phi$ can be extended on monomials as the axioms require. The fact that $(PR+1)(PQ-1)=0$ implies that $\Delta$ can be extended as an algebra homomorphism. Finally, the fact that $PR$ commutes with $PQ$ imples that the antipode $S$ defined on generators can be extended using $S\,m=m(S\otimes S)\Phi$ and
$\Phi(S\otimes \id)=(\id\otimes S)\Phi$. This and similar results are also to be found in [HH]. $\square$\smallbreak

\noindent  We observe that the relations involving $Q_{ij}$ and $R_{kl}$ ensure that $\cov$ and $\vec$ are comodule algebras for $A(R)$. Hence if $A(R)$ can be extended or quotiented to a dual quasitriangular Hopf algebra $H(R)$ with  dual quasitriangular structure consistent with that of $A(R)$, and if $\cov$ and $\vec$ are comodule algebras for $H(R)$, then all possible vector spaces obtained by tensoring $\cov$ and $\vec$ can be provided of an algebra structure such that they are again comodule algebras for $H(R)$. The product is then defined by means of the braidings and the product in $\cov$ and $\vec$.\smallbreak

\noindent We want to study the case in which $R$ is (a multiple of) the matrix solution of the quantum Yang-Baxter equation associated to type $A_n$ (see [FRT]) with real entries because of our choice of $q$. We regard the braided covector and vector algebras as complex algebras. Since $R$ is Hecke, after some rescaling we have:

\noindent $R^{jj}_{jj}=q^2$, $R^{ij}_{ij}=q$ if $i\not=j$, $R^{ij}_{ji}=q^2-1$ if $i<j$,
and all the other entries equal to $0$.

\noindent We choose $Q=q^{-2}R$, i.e.
$Q^{jj}_{jj}=1$, $Q^{ij}_{ij}=q^{-1}$ if $i\not=j$, $R^{ij}_{ji}=1-q^{-2}$ if $i<j$, and all the other entries equal to $0$.  It is well known that in this case the matrix $R$ is bi-invertible, i.e. $R$ and $R^{t_2}$ are both invertible, where  $t_2$ stands for transposing only with respect to the second component in the tensor product, i.e. there exists an $n^2\times n^2$ matrix
$\tilde R=((R^{t_2})^{-1})^{t_2}$ such that $(\tilde R)^{aj}_{ib}R^{id}_{cj}=R^{aj}_{ib}(\tilde R)^{id}_{cj}=\delta^a_c\delta_b^d$.

\noindent Then the relations and the braiding for $\cov$ are
$x_ix_j=qx_jx_i$ if $i>j$ so this is the $n-$dimensional quantum space. \noindent We know that this algebra has a basis given by (increasing) ordered monomials $x_1^{e_1}\cdots\,x_n^{e_n}$. The braiding is given on the generators by
$$\Phi(x_i\otimes x_j)=\cases{qx_j\otimes x_i&if $i<j$\cr
q^2x_i\otimes x_i& if $i=j$\cr
(q^2-1)x_i\otimes x_j+qx_j\otimes x_i& if $i>j$\cr}$$
so that, for $i<j$, it follows that $\Phi(x_i^a\otimes x_j^b)=q^{ab}x_j^b\otimes x_i^a$. We use this, the properties and the $q$-binomial formula to see that
the comultiplication of an arbitrary element of the basis equals:
$$\Delta(x_1^{e_1}\cdots\,x_n^{e_n})=\sum_{j_1=0}^{e_1}\,\cdots\,\sum_{j_n=0}^{e_n}\left(\prod_{j=1}^n\left[e_i\atop j_i\right]_{q^2}\right)
q^{\sum_vJ_v(e_v-j_v)}x_1^{e_1-j_1}\cdots\,x_n^{e_n-j_n}\otimes x_1^{j_1}\cdots\,x_n^{j_n}$$
\noindent The antipode is given by $S(x_1^{e_1}\cdots\,x_n^{e_n})=q^{|E|^2-|E|}
(-1)^{|E|}x_1^{e_1}\cdots\,x_n^{e_n}$ as one can easily prove by induction.\vskip0.5cm

\noindent For convenience we put  $\partial_j:=\partial^j$. The relations for $\vec$ are given by $\partial_i\partial_j=q\partial_j\partial_i$ if $i<j$. The braiding is given on generators by
 
$$\Phi(\partial_i\otimes \partial_j)=\cases{q\partial_j\otimes \partial_i&if $i>j$\cr
q^2\partial_i\otimes \partial_i& if $i=j$\cr
(q^2-1)\partial_i\otimes \partial_j+q\partial_j\otimes \partial_i& if $i<j$\cr}$$ and can be extended on monomials according to the axioms. Again we know that ordered monomials provide a basis for $\vec$.
In this case we fix the basis given by ordered monomials where the order is {\sl decreasing}. 

\noindent One can observe that for every choice of nonzero constants $c_j$, for $j=1,\ldots,\,n$,  there is an algebra isomorphism $\psi$ between $\cov$ and $\vec$ mapping $x_j$ to $c_j\partial_{n+1-j}$, such that $\Delta_{\vec}\circ\psi=(\psi\otimes\psi)\circ\Delta_{\cov}$ and $S_{\vec}\circ\psi=S_{\cov}$. In particular, for $c_j=q^{-j+{1\over2}(n-1)}$ for every $j$, then we also have $\Phi_{\vec,\vec}(\psi\otimes\psi)=\Phi_{\cov,\cov}$. 

\noindent Using this property, or by computations similar as for $\cov$, one finds the comultiplication and the antipode on basis elements: 
$$\Delta(\partial_n^{e_n}\cdots\,\partial_1^{e_1})=\sum_{j_n=0}^{e_n}\,\cdots\,\sum_{j_1=0}^{e_1}\left(\prod_{j=1}^n\left[e_i\atop j_i\right]_{q^2}\right)
q^{\sum_uJ^u(e_u-j_u)}\partial_n^{e_n-j_n}\cdots\,\partial_1^{e_1-j_1}\otimes \partial_n^{j_n}\cdots\,\partial_1^{j_1}$$
$$S(\partial_n^{e_n}\cdots\,\partial_1^{e_1})=q^{|E|^2-|E|}(-1)^{|E|}\partial_n^{e_n}\cdots\,\partial_1^{e_1}$$\vskip0.5cm
\noindent $A(R)$ can either be quotiented to the algebra of functions on quantum $SL_n$ or localized at the quantum determinant in order to obtain the algebra of functions on quantum $GL_n$, and both are Hopf algebras. Hence by  Majid's theory (see Corollary 9.2.14 and Proposition 10.3.6 in [Ma]) we recover the braiding between $\cov$ and $\vec$ and between $\vec$ and $\cov$.
They are given by $\Phi(\partial_i\otimes x_j)=\sum_{r,s}x_r\otimes\partial_s(R^{-1})^{ir}_{sj}$ and $\Phi(x_i\otimes\partial_j)=\sum_{r,s}(\tilde R)^{rj}_{is}\partial_s\otimes x_r$. By the simplicity of the case $A_n$, one can also check directly that these define braidings such that the two  braided tensor products are again  braided algebras (i.e. algebras with a braiding with the necessary compatibility). 
 We compute $R^{-1}$ and $\tilde R$ and we see that:
$$\Phi_{\vec,\cov}(\partial_i\otimes x_j)=\cases{q^{-1}x_j\otimes \partial_i&if $i\not=j$\cr
q^{-2}x_j\otimes \partial_j+\sum_{r>j}(q^{-2}-1)x_r\otimes\partial_r& if $i=j$\cr}$$ and

$$\Phi_{\cov,\vec}(x_i\otimes \partial_j)=\cases{q^{-1}\partial_j\otimes x_i&if $i\not=j$\cr
\sum_{r<j}(q^{-2}-1)q^{-2(j-r)}\partial_r\otimes x_r+q^{-2}\partial_j\otimes x_j& if $i=j$\cr}$$
\noindent We observe that for the algebra isomorphism $\psi$ described above, with $c_j=q^{-j+{1\over2}(n-1)}$ we have:
$$(\psi^{-1}\otimes\psi)\circ\Phi_{\cov,\vec}=\Phi_{\vec,\cov}\circ(\psi\otimes\psi^{-1})$$ and $$(\psi\otimes\psi^{-1})\circ\Phi_{\vec,\cov}=\Phi_{\cov,\vec}\circ(\psi^{-1}\otimes\psi)$$ However, $\psi$ cannot be considered as a morphism in the braided category since it is {\sl not true} that $(\psi\otimes \id)\circ\Phi_{\cov,\vec}=\Phi_{\cov,\cov}\circ(\id\otimes\psi)$, as one can easily see by computing lthe actions of the left hand side and of the right hand side on $(\partial_2\otimes x_1)$ for $n=2$.\vskip0.2cm

\noindent Majid introduces also a theory of differentiation. Namely, every $\partial_j$ acts
on $\cov$ as braided partial differentiation. Majid defines this action as a sort of analogue of the  linear approximation of a function. Namely, for a basis element $\xe$, one computes $\Delta(\xe)\in\cov\otimes\cov$ and one expresses the latter as sum of elements as ``basis elements of $\cov$ tensor elements of $\cov$''. Then $\partial_i(\xe)$ is the term appearing as a ``coefficient'' of $x_i\otimes 1$. One can easily check that this defines an action of $\vec$ on $\cov$ and one can compute that for every $f(\ux)\in\cov$  this turns out to be:
$$\partial_jf(\ux)=x_j^{-1}\left[{f(q^2x_1,\,\ldots,\,q^2x_j,\,x_{j+1},\ldots,\,x_n)-f(q^2x_1,\,\ldots,\,q^2x_{j-1},\,x_j,\,x_{j+1},\ldots\,x_n)\over(q^2-1)}
\right]$$ where
the inverse of $x_j$ is only formal, and ``apparent''. However, if the reader doesn't like
this, he can work in the Ore localization of $\cov$ at $x_j$, since
the set $\{x_j^k\;|\;k\ge0\,\}$ satisfies Ore conditions.\smallbreak
\noindent In particular, one has $\partial_j(\xe)=[e_j]_{q^2}q^{E_j}x_1^{e_1}\cdots x_{j-1}^{e_{j-1}}x_j^{e_j-1}x_{j+1}^{e_{j+1}}\cdots x_n^{e_n}$.\smallbreak 
\noindent Formally we can work also with $\cove$ (respectively $\vece$), the algebra of formal power series in the $x_j$'s (respectively $\partial_j$'s) with the given defining relations. In this case everything that we have described above  is defined as in $\cov$ and $\vec$.

\beginsection3. The exponential map, $e_q$, and $E_q$

\noindent As in [Ma], we introduce the exponential map as a coevaluation for $\cov$ and $\vec$. We say in a few lines how this works. 

\noindent There is a braided Hopf algebra pairing $<\,,\,>\colon \cov\otimes\vec\to K$ where $K$ denotes the field where the braided Hopf algebras are defined, given by $$<g(\partial_n,\,\ldots,\,\partial_1),\,f(x_1,\,\ldots,\,x_n)>=\varepsilon(g(\partial_n,\,\ldots,\,\partial_1)f(x_1,\,\ldots,\,x_n))$$ (see [Ma] for further details). The dual basis to the one given by $x_1^{e_1}\cdots x_n^{e_n}$ is then given by $(\xe)^*:={{\partial_n^{e_n}}\over {[e_n]_{q^2}!}}\cdots\,{{\partial_1^{e_1}}\over {[e_1]_{q^2}!}}$. The coevaluation map $\exp\colon K\to(\cov\otimes\vec)^{\rm ext}$ is uniquely determined by the image $\exp(x|\,\partial)$ of $1$ in a formal extension to power series of $\cov\otimes\vec$. The element $\exp(x|\,\partial)$  is defined as the canonical element $\sum_{E,F}c_{EF}\xe\otimes\df$ such that for every $f(\ux)$ in $\cov$ and every $g(\upar)$ in $\vec$, one has:
$$\eqalign{&f(\ux)=\sum_{E,F}c_{EF}\xe<\df,\,f(\ux)>\quad{\hbox{ and}}\cr
& g(\upar)=\sum_{E,F}c_{EF}<g(\upar),\,\xe>\df\cr}$$ Hence it is equal to
$$\sum_{e_1,\,\ldots,\,e_n\ge0}x_1^{e_1}\cdots x_n^{e_n}\otimes(\xe)^*=\sum_{e_1,\,\ldots,\,e_n\ge0}x_1^{e_1}\cdots x_n^{e_n}\otimes{{\partial_n^{e_n}}\over {[e_n]_{q^2}!}}\cdots\,{{\partial_1^{e_1}}\over {[e_1]_{q^2}!}}$$ 
\noindent By Example 10.4.16 in [Ma] this is equal to $e_{q^{-2}}((1-q^{-2})\sum_{i=1}^nx_i\otimes\partial_i)$ where $e_{q}(z)=\sum_{k=0}^{\infty}{{z^k}\over{(q;q)_k}}$  (see [Koo] for further details). It follows by straightforward computation that $\exp(x\,|\,\partial)$ is also equal to $E_{q^2}((1-q^{2})\sum_{i=1}^nx_i\otimes\partial_i)$, where $E_{q^2}(z)=\sum_{k=0}^{\infty}{{q^{{1\over2}k(k-1)}z^k}\over{(q;\,q)_k}}$. 
Properties analogous to those of the classical exponential map hold for $\exp$ as a consequence of the fact that it is a canonical element arising from a braided Hopf algebra pairing. In our case the properties we need can be also directly checked knowing the properties of $e_{q}(z)$ and $E_q(z)$. Corollary 10.4.17 in [Ma] gives us also a braided version of Taylor formula, which can be constructed also easily directly in the $A_n$ case. This is given by $$\eqalign{\Delta(f(x_1,\ldots,\,x_n))&\cr
&=f(\Delta(x_1),\ldots,\,\Delta(x_n))=f(x_1+y_1,\,\ldots,\,x_n+y_n)\cr
&=\sum_{e_1,\,\ldots,\,e_n\ge0}x_1^{e_1}\cdots x_n^{e_n}\otimes\biggl({{\partial_n^{e_n}}\over {[e_n]_{q^2}!}}\cdots\,{{\partial_1^{e_1}}\over {[e_1]_{q^2}!}}f(x_1,\ldots,\,x_n)\biggr)\cr
&=\sum_{e_1,\,\ldots,\,e_n\ge0}y_1^{e_1}\cdots y_n^{e_n}\biggl({{\partial_n^{e_n}}\over {[e_n]_{q^2}!}}\cdots\,{{\partial_1^{e_1}}\over {[e_1]_{q^2}!}}f(x_1,\ldots,\,x_n)\biggr)\cr
&=\exp(y\,|\,\partial)f(\ux)}$$ where $y_j=x_j\otimes 1$ and $x_i$ stands for $1\otimes x_i$ after the second equality sign. In the first equality we made an abuse of notation identifying the power series $f$ with the series of its coefficients.

\beginsection4. Integration on $\cove$

In [KeMa] a theory of integration on braided covector and braided vector Hopf algebras is introduced. The authors work out a few particular cases in more detail. Integration on the braided line was futher worked out in [Koo] who  proved, by means of finding  suitable representations, the left invariance of the integral in an analytical way. The main tool was finding a representation for
(an extension of ) the algebra of ``quantum coordinates'' on the braided line, and a common eigenvector for the action. In his case, an eigenvector for the integral of an element of the braided line would have a Jackson integral in commuting variables
as eigenvalue. Therefore, convergence and equalities were  investigated on the eigenvalues. This is no longer possible if the dimension of the braided space is bigger than one, because noncommuting operators cannot have a common eigenvector. We have to work a little bit more to give a meaning to our version of the integral. In many cases we see that the concept of integrability may depend on the choice of the representation we choose. This is equivalent to the choice of a normal form in our case.

\noindent We start with the ``indefinite'' integral with respect to $x_i$.
As in [KeMa], we view it as an operator on $\cov$.

\noindent \proclaim Definition. The braided partial integral with respect to $x_i$ acting on $f(x_1,\ldots\,x_n)\in\cove$ is given by:
$$\int^{x_i}_0
f:=(1-q^2)\sum_{k=0}^\infty q^{2k}x_if(q^{-2}x_1,\,\ldots\,q^{-2}x_{i-1},\,q^{2k}x_i,\,x_{i+1},\ldots\,x_n)$$

\noindent It's easy to see that the operator defined above acts as a pseudo
inverse for the partial differential operator $\partial_i$. It is indeed only a right inverse, since it acts as a left inverse for $\partial_i$ only on series containing $x_i$ in every monomial of its expansion (see remark in [KeMa], page 6815). Each $\int_0^{x_i}$ is a well defined operator from $\cove$ to $\cove$ since for every basis monomial $\xe$ we can write $\int_0^{x_i}\xe$ as a monomial in the $x_j$'s with a coefficient that is a {\sl convergent} series of complex numbers for $0<q<1$.

\noindent Following [KeMa] one could write this operator as a power series, namely:
$$(1-q^2)x_i\sum_{k=0}^{\infty}(1-(1-q^2)\partial_ix_i)^k$$ which makes sense because the usual basis is a basis of eigenvectors, with convergent eigenvalues.

\noindent We can read $\int^{x_i}_0f$ as a ``function'' of the $x_j$'s, hence it makes sense  to consider, for instance $\int_0^{ax_i}$ for a nonzero constant $a$. In particular we can define $\int^{-x_i}_0f$ and $\int^{q^lx_i}_0f$ as follows.

$$\int^{-x_i}_0
f:=-(1-q^2)\sum_{k=0}^\infty q^{2k}x_if(q^{-2}x_1,\,\ldots\,q^{-2}x_{i-1},\,-q^{2k}x_i,\,x_{i+1},\ldots\,x_n)$$
$$\int^{q^lx_i}_0
f:=(1-q^2)\sum_{k=0}^\infty q^{2k+l}x_if(q^{-2}x_1,\,\ldots\,q^{-2}x_{i-1},\,q^{2k+l}x_i,\,x_{i+1},\ldots\,x_n)$$

\noindent We define then  for every $f\in\cove$
$$\int_{-x_i}^{x_i}f:=\int_0^{x_i}f-\int_0^{-x_i}f$$ and
$$\int_{-x_i\cdot\infty}^{x_i\cdot\infty}f:=\lim_{r\to\infty}\int_{-q^{-2r}x_i}^{q^{-2r}x_i}f$$
\noindent The last definition is only formal so far, because the image of a power series is no longer a power series (coefficients might be infinite sums themselves), and we have no notion of convergence. Indeed we cannot find a convergence set because we cannot give values to the variables, since they do not commute. This issue can be solved in different ways, so that we can give a meaning to equalities as well. The ideas here are based on the apporaches of Kempf and Majid, and of Koornwinder. We approach the problem not for a single infinite integral, but for the $n-$dimensional integral:
$$I\cdot f:=\int_{-x_n\cdot\infty}^{x_n\cdot\infty}\cdots \,\int_{-x_1\cdot\infty}^{x_1\cdot\infty}f=$$ formally $$(1-q^2)^n\sum_{k_n=-\infty}^{\infty}\cdots\,\sum_{k_1=-\infty}^{\infty}\sum_{{\bf\varepsilon}\in\{\pm1\}^n}q^{2|K|+{n\choose2}}x_1\cdots\,x_n\,f\,(\varepsilon_1q^{2k_1}x_1,\,\ldots,\,\varepsilon_nq^{2k_n}x_n)$$ \noindent As we said $I\cdot f$ is no longer an element of $\cov$, since the coefficients with respect to the elements of the basis are not always definite. We fix then an action of $\cov$ on the space of power series in the $n$ commuting variables $z_1,\,\ldots,\,z_n$ with complex coefficients. This representation corresponds with the choice of a normal form for the monomials in $\cov$. The representation, denoted by $\triangleright$, for monomials in the $x_i$'s acting on monomials in the $z_j$'s, is given by:
$$\xe\triangleright z_1^{h_1}\cdots\,z_n^{h_n}=q^{\sum_{i=1}^ne_iH_i}z_1^{h_1+e_1}\cdots\,z_n^{h_n+e_n}$$ and can be extended linearly to an action of $\cov$ on formal power series in the $z_j$'s. We can restrict the space on which we act by taking the space $V$ of power series which are absolutely convergent in a neighbourhood of zero. This makes sense because the $z_i$'s commute with each other. We see that this space is invariant under the action of $\cov$. Moreover we see that we can extend the representation of $\cov$ on $V$ to a representation of the class $C$ given by the power series $f$ in the $x_i$'s such that $f\triangleright 1\in V$. Indeed one can see that

\noindent (1) $\forall\,f=f(x_1,\,\ldots,\,x_n)\in C$ and $\forall\, g=g(z_1,\,\ldots,\,z_n)\in V$ it holds that $(f\triangleright g)(z_1,\,\ldots,\,z_n)$ belongs to $V$ because the associated  series of absolute values is majorized by the product in $V$ of series of absolute values associated to $f\triangleright 1$ and $g$.

\noindent (2) $\forall\,f$ and $g\in C$, their product $fg\in C$ because $(fg)\triangleright 1=f\triangleright(g\triangleright1)\in V$ by (1)

\noindent Moreover,

\noindent (3) $\forall f\in C$, $\quad f\triangleright 1=0$ $\Leftrightarrow$ $f\equiv0\quad$ i.e. all the coefficients of $f$ have to be zero. 

\noindent From now on we write $f_{\cdot1}$ for $f\triangleright 1$, for any expression $f(\ux)$ for which the action on $1$ makes sense.

\noindent We would like to extend the representation now to the formal expressions of type $I\cdot f$ for $f\,\in C$. This does not always make sense, hence we have to add further conditions. Let us take the subclass $C'$ of $C$ given by the series in $C$ such that 

\noindent (a) $f_{\cdot1}$ can be continued analytically on ${\bf R}^n+iU$ for some open neighboorhood $U$ of $0$ in ${\bf R}^n$;  

\noindent (b) $f_{\cdot1}$ is absolutely $q^2$-integrable for every $\uz\in{\bf R}^n$ for which every $z_j\not=0$, i.e.
$$(1-q^2)^n\sum_{k_1=-\infty}^{\infty}\cdots\,\sum_{k_n=-\infty}^{\infty}\,\sum_{{\bf \varepsilon}\in \{\pm1\}^n}q^{2|K|}|z_1|\cdots\,|z_n|\,|f\cuno(\varepsilon_1q^{2k_1}z_1,\ldots,\,\varepsilon_nq^{2k_n}z_n)|<\infty$$ for $\uz$ outside the standard hyperplanes.

\noindent Then we can say, for $f\in C'$ what is $(I\cdot f)\cuno$:

$$\eqalign{&{{(I\cdot f)\cuno}}=\cr
&=\left[\int_{-x_n\cdot\infty}^{x_n\cdot\infty}\cdots \,\int_{-x_1\cdot\infty}^{x_1\cdot\infty}f\right]\triangleright 1\cr
&=\left[(1-q^2)^n\sum_{k_n=-\infty}^{\infty}\cdots\sum_{k_1=-\infty}^{\infty}\sum_{{\bf\varepsilon}\in\{\pm1\}^n}q^{2K+{n\choose2}}x_1\cdots\,x_n\,f\,(\varepsilon_1q^{2k_1}x_1,\ldots,\varepsilon_nq^{2k_n}x_n)\right]\triangleright 1\cr}
$$

\noindent By expressing every term in the $x_i$'s as a sum of elements of the basis (with in principle infinite sums as coefficients) and applying the action, we see that the coefficients make sense and that

$$\eqalign{&{(I\cdot f)\cuno}=\cr
&=q^{n\choose 2}\int_{-z_1\cdot\infty}^{z_1\cdot\infty}\cdots \,\int_{-z_n\cdot\infty}^{z_n\cdot\infty}f\cuno(q^{n-1}t_1,\,\ldots,q^{(n-i)}t_i,\,\ldots,\,t_n)d_{q^2}t_n\cdots
d_{q^2}t_1=\cr
&=\int_{-q^{n-1}z_1\cdot\infty}^{q^{n-1}z_1\cdot\infty}\cdots \int_{-q^{n-i}z_i\cdot\infty}^{q^{n-i}z_i\cdot\infty}\cdots\,\int_{-z_n\cdot\infty}^{z_n\cdot\infty} f\cuno(t_1,\,\ldots,\,t_n)d_{q^2}t_n\cdots
d_{q^2}t_1\cr}$$ where the $q^2$-integral in the second line is the $q^2$-Jackson integral in $n$ variables obtained by iterating $(8.11)$ in [Koo].
It is clear that if $f\in C'$ then $(I\cdot f)\cuno$ converges whenever $z_j\not=0$ for every $j$.

\noindent Because of (3),  two objects in $C$ are equal if and only if they act in the same way on $1$. Following this philosophy, we say that two   $q^2$-integrals of objects in $C'$ are equal if and only if they act on the same way on $1$. This will be our tool to show equalities then. 

\noindent The first purpose is to show in a less formal way than in [KeMa] where this appeared first, the invariance of the operator $I$ under translation. As in [KeMa] and in [Koo], we make use of the Taylor formula. For this we need an extra assumption on the elements in $C'$, since we have to use partial derivatives. We consider $f\in C'$ satisfying:

\noindent (c) For some $\eta>0$ there exists for each $J\in({\bf Z}_{\ge0})^n$ some constant $C_J$ such that
$$|(D_{1,q^2}^{j_1}\cdots D_{n,q^2}^{j_n}f\cuno)(z_1,\ldots,\,z_n)|\le C_J\prod_{k=1}^n(1+|z_k|^2)^{-(1+\eta)}$$ if $\uz\in{\bf R}^n$.

\noindent For an $f$ in $C'$, condition (c) implies that all the  Jackson  derivatives of $f\cuno$ are absolutely $q^2$-integrable for all $\uz$ in a neighbourhood of $0$ minus the intersection with the standard hyperplanes.

\noindent One sees immediately that $(\partial_if)\cuno=(D_{i,q^2} f\cuno)(qz_1,\ldots,\,qz_{i-1},\,z_i\ldots,\,z_n)$.

\noindent We show that it makes sense to compute $I(\partial_n^{j_n}\cdots\partial_1^{j_1}f)\triangleright 1$, and that this is equal to zero whenever $(j_1,\,\ldots,\,j_n)\not=(0,\ldots,\,0)$.
 
\noindent We write:
 \noindent$F^J\cuno(z_1,\ldots,\,z_n):=(\partial_n^{j_n}\cdots\partial_1^{j_1}f)\triangleright 1=(D_{1,q^2}^{j_1}\cdots\,D_{n,q^2}^{j_n} (f\cuno))(q^{J^1}z_1,\ldots,\,q^{J^n}z_n)$. One sees that $F^J\cuno\in V$ if $f$ satisfies condition (c). As far as $q^2$-integrability is concerned: $F^J\cuno$ is absolutely $q^2$- Jackson integrable if and only if for every choice of $(h_1,\ldots,\,h_n)\in \{\pm1\}^n$ $$\sum_{{\underline\varepsilon}\in\{\pm1\}^n}\sum_{k_1=0}^{\infty}\cdots\sum_{k_n=0}^{\infty}q^{2K\cdot H}|z_1\cdots\,z_n|\,|F^J\cuno(q^{2k_1h_1}\varepsilon_1z_1,\ldots,\,q^{2k_nh_n}\varepsilon_nz_n)|$$ has a positive radius of convergence. Hence, if condition (c) holds for $f(\ux)$, and since the above sums are of the form:
$$\eqalign{&\sum_{\underline\varepsilon}\sum_{k_1=0}^{\infty}\cdots\sum_{k_n=0}^{\infty}q^{2K\cdot H}|z_1\cdots z_n\,(D_{1,q^2}^{j_1}\cdots D_{n,q^2}^{j_n}f\cuno)
(q^{2k_1h_1+J^1}\varepsilon_1z_1,\ldots,q^{2k_nh_n+J^n}\varepsilon_nz_n)|\cr
&\le C_J|z_1\cdots z_n\,|\sum_{\underline\varepsilon}\sum_{k_1=0}^{\infty}\cdots\sum_{k_n=0}^{\infty}q^{2K\cdot H}\prod_{r=1}^n(1+|q^{2k_rh_r}z_r|^2)^{-(1+\eta)}\cr}$$ that converges since $0<q<1$, one sees that the $F\cuno^J$ are $q^2$-integrable. 

\noindent Moreover for $\gamma=(\gamma_1,\,\ldots,\,\gamma_n)$ with $\gamma_j$ in ${\bf R}-\{0\}$ we have
$$\eqalign{&q^{\sum_{k=1}^nJ^k}\int_{-\gamma_1\cdot\infty}^{\gamma_1\cdot\infty}\cdots\int_{-\gamma_n\cdot\infty}^{\gamma_n\cdot\infty}|F^J\cuno(t_1,\ldots,\,t_n)|d_{q^2}t_n\cdots d_{q^2}t_1\cr
&=\int_{-q^{J^1}\gamma_1\cdot\infty}^{q^{J^1}\gamma_1\cdot\infty}\cdots\int_{-q^{J^n}\gamma_n\cdot\infty}^{q^{J^n}\gamma_n\cdot\infty}|D_{1,q^2}^{j_1}\cdots\,D_{n,q^2}^{j_n}f\cuno(t_1,\ldots,\,t_n)|d_{q^2}t_n\cdots d_{q^2}t_1=0\cr}$$ The proof is as in the one-dimensional case (see [Koo]). 

\noindent Hence we can state that
\proclaim Lemma 4-1. Let $f\in C'$ satisfy condition {\rm (c)}. Then for every $J\not=\underline0$ there holds: $(I(\partial_n^{j_n}\cdots\partial_1^{j_1}f))\triangleright 1\equiv0$, so that we can conclude that $I(\partial_n^{j_n}\cdots\partial_1^{j_1}f)=0$ .\par

\noindent{\bf Proof: } One has
$$\eqalign{&(I(\partial_n^{j_n}\cdots\partial_1^{j_1}f))\cuno=\cr
&=\int_{-q^{n-1}z_1\cdot\infty}^{q^{n-1}z_1\cdot\infty}\cdots \int_{-q^{n-i}z_i\cdot\infty}^{q^{n-i}z_i\cdot\infty}\cdots\,\int_{-z_n\cdot\infty}^{z_n\cdot\infty}F^J\cuno(t_1,\,\ldots,\,t_n)d_{q^2}t_n\cdots
d_{q^2}t_1=0\qquad\square\cr}$$

\proclaim Proposition 4-2. Let $f\in C'$ satisfy {\rm (c)}. Then $(id\otimes I)\Delta(f)=1\otimes (If)$.\par
\noindent{\bf Proof: } By the braided Taylor formula we have
$$(id\otimes I)(\Delta f)=\sum_{j_1,\ldots,\,j_n\ge0}{{y_1^{j_1}\cdots\,y_n^{j_n}}\over{[j_1]_{q^2}!\,\cdots[j_n]_{q^2}!}}I(\partial_n^{j_n}\cdots\partial_1^{j_1}f)$$ but this is by Lemma 4-1 equal to the term with all $j_k$'s equal to $0$. $\square$.

\proclaim Proposition 4-3. Let $f(\ux)\in C$ such that {\rm (c)} holds for every $\eta>0$. Then the statement of Proposition 4-2 holds for every polynomial $p(\ux)$ times $f(\ux)$.\par
\noindent{\bf Proof: } It is not restrictive to assume that $p(\ux)$ is a monomial. By (1) it follows that also $p(\ux)\cdot f(\ux)\in C$. We have to check that condition (c) holds for every element of the form $\xe f({\underline x})$. One sees immediately that 
$$(\xe f({\underline x}))\triangleright 1=f\cuno(q^{E^1}z_1,\,\ldots,\,q^{E^n}z_n)z_1^{e_1}\cdots z_n^{e_n}\in V$$
if $f(\ux)\in C$, so that $(\xe\,f(\ux))\cuno$ makes sense. 

\noindent Condition (c) is on the $q^2$-Jackson partial derivatives on commuting variables, and it holds for $\xe f({\ux})$ as a consequence of the fact that for two functions $a({\uz})$ and $b({\uz})$ and for any $j=1,\ldots,\,n$:
$$D_{j,q^2}(a({\uz})b({\uz}))=(D_{j,q^2}a({\uz}))(b(z_1,\ldots,\,z_{j-1},\,q^2z_j,\,z_{j+1},\,\ldots,\,z_n))+a({\uz})D_{j,q^2}(b({\uz}))$$ Then, by Proposition 4-2 we have the statement. $\square$ 

\noindent We have a description of a class of power series in the $x_i$'s for which integration makes sense, although we are not able so far to make a complete classification of integrable functions. The same problem is treated in [OR] for the one-dimensional case. Still, what we have is enough to allow computations in the following case.\smallbreak

\noindent{\bf Example 1:} The $q^2$-Gaussian $g_{q^2}(\ux)$ is $$g_{q^2}(\ux):=e_{q^4}(-\ux\cdot\ux)=e_{q^4}(-\sum_{j=1}^nx_j^2)=\prod_{j=1}^ne_{q^4}(-x_j^2)$$ where the last equality holds because of Proposition 3.1 in [Koo] (this  result is due to [Sch]) and the product in the above formula is taken with increasing order on the variables. It satisfies conditions (a), (b) and (c) for every $\eta>0$. This is a consequence of the one dimensional case (see [Koo]) and the fact that $$(g_{q^2}(\ux))\cuno(\uz)=\prod_{j=1}^ne_{q^4}(-z_j^2)\qquad\partial_j(g_{q^2}(\ux))=-{{x_j}\over{(1-q^2)}}g_{q^2}(\ux)$$ and that $$D_{q^2,1}^{j_1}\cdots\,D_{q^2,n}^{j_n}(e_{q^4}(-\sum_{k=1}^n z_k^2))=p(\uz)\,e_{q^4}(-\sum_{k=1}^n z_k^2)$$ where $p(\uz)$ is a polynomial in the $z_j$'s. 

\noindent It follows then that also elements of the form: $$x_1^{a_1}g_{q^2}(x_1)\cdots x_n^{a_n}g_{q^2}(x_n)=x_1^{a_1}\cdots x_n^{a_n}e_{q^4}(-\sum_j(q^{-A^j}x_j)^2)$$ satisfy condition (c), so that for every $p_j(x_j)\in\cov$ we can integrate every element of the form $p_1(x_1)g_{q^2}(x_1)\cdots p_n(x_n)g_{q^2}(x_n)$.\smallbreak
\noindent In particular, for $f(\ux)=g_{q^2}(x_1)x_1^{a_1}\cdots g_{q^2}(x_n)x_n^{a_n}$
one can compute $$\eqalign{&(I\cdot f)\cuno|_{\uz=\gamma}\cr
&=\int_{-q^{n-1}\gamma_1\cdot\infty}^{q^{n-1}\gamma_1\cdot\infty}\cdots \int_{-q^{n-i}\gamma_i\cdot\infty}^{q^{n-i}\gamma_i\cdot\infty}\cdots\,\int_{-\gamma_n\cdot\infty}^{\gamma_n\cdot\infty} f\cuno(t_1,\,\ldots,\,t_n)d_{q^2}t_n\cdots
d_{q^2}t_1\cr
&=\prod_{j=1}^n\left(\int_{-q^{n-j}\gamma_j\cdot\infty}^{q^{n-j}\gamma_j\cdot\infty}e_{q^4}(-t_j^2)t_j^{a_j}d_{q^2}t_j\right)\cr
&=\cases{\prod_{j=1}^n\biggl(c_{q^2}(\gamma_jq^{n-j})q^{-{{a_j^2}\over 2}}(q^2;q^4)_{{a_j}\over2}\biggr)& if $a_j$ even $\forall j$\cr
0& otherwise\cr}\cr}$$ where $$c_{q^2}(\gamma)={{2(1-q^2)(q^4,\,-q^2\gamma^2,\,-q^2\gamma^{-2};\,q^4)_{\infty}}\over{(-\gamma^2,\,-q^4\gamma^{-2},\,q^2;\,q^4)_{\infty}}}$$ as in formula $(8.15)$ in [Koo]. In particular 
$$\bigl(I\cdot f\bigr)\cuno|_{\uz=\gamma}=q^{-\sum_ja_j^2/2}\prod_j(q^2;q^4)_{a_j/2}\bigl(I\cdot g_{q^2}(\ux)\bigr)|_{\uz=\gamma}$$ if all $a_j$'s are even, and $0$ otherwise, hence we can conclude that:$$I\cdot f=\cases{q^{-\sum_ja_j^2/2}\prod_j(q^2;q^4)_{a_j/2}\bigl(I\cdot g_{q^2}(\ux)\bigr)& if $a_j$ is even $\forall j$\cr
0& otherwise\cr}$$
The observation that the integral of the Gaussian times a monomial is equal to a constant times the integral of the Gaussian appeared in [KeMa] and in [Ma] first. In our case we have to deal with the shift and this depends on the choice of our global integral.$\;\spadesuit$\smallbreak

\beginsection5. Lattice integrability

\noindent In the previous Section we saw a definition of integrable series in $\cove$. Unfortunately the above method fails for another analogue
of the Gaussian we would like to $q^2$-integrate, namely the $q^2$-Gaussian $$G_{q^2}(\ux):=E_{q^4}(-\ux\cdot\ux)=E_{q^4}(-\sum_{j=1}^nx_j^2)=E_{q^4}(-x_n^2)\cdots E_{q^4}(-x_1^2)$$ where the last equality holds because of the result in [Sch] to be found in [Koo], Proposition 3.1. In [Koo], section 9, it is also shown  that $G_{q^2}(\ux)$ does not satisfy condition (c), nor condition (b) even for $n=1$. On the other hand it is also shown there that for a given choice of a $q^2$-lattice of the form $\{q^{2k}\gamma\,|\;k\in{\bf Z}\}$, namely for $\gamma=1$,  $(I\cdot G_{q^{2}}(x_1))\cuno\;|_{z_1=1}$ is absolutely convergent. Hence one can introduce a weaker version of integrability in $\cove$, that we will call {\sl ``lattice integrability''} requiring for an $f(\ux)$  such that $f\cuno$ is entire, that there is a $q^2$-lattice $L(\gamma)=\{\,(\gamma_1q^{2k_1},\ldots,\,\gamma_nq^{2k_n})\ |\; K\in{\bf Z}^n\}$ in  ${\bf R}_{\not=0}^n$ such that $(I\cdot f)\cuno\,|_{\uz\in L(\gamma)}$ is absolutely convergent. One can easily see that the series $f(\ux)=E_{q^4}(-x_1^2)\cdots E_{q^4}(-x_n^2)$ is lattice-integrable for $\gamma=(q^{n-1},\,\ldots,\,q^{n-j},\,\ldots,1)$. This $f$ is unfortunately not the $q^2$-Gaussian $G_{q^2}(\ux)$ that we wanted to integrate, for $n\ge2$. Moreover, we can show that already for $n=2$, $G_{q^2}(\ux)$ is not lattice integrable although it is entire. \smallbreak

\noindent{\bf (Counter)example 1:} Let us consider $G_{q^2}(\ux)$ for $n=2$.
We write for simplicity $x_1=x$ and $x_2=y$, and $z_1=z$, $z_2=w$.
Then $$G_{q^2}(\ux)=\sum_{k,l=0}^\infty {{(-1)^{k+l}q^{2k^2+2l^2-2k-2l}y^{2l}x^{2k}}\over{(q^4;q^4)_l\,(q^4;q^4)_k}}=\sum_{l=0}^\infty{{(-1)^{l}q^{2l^2-2l}E_{q^4}(-(q^{2l}x)^2)y^{2l}}\over{(q^4;q^4)_l}}$$ hence $$\eqalign{&G_{q^2}(\ux)\cuno=\sum_{k,l=0}^\infty {{(-1)^{k+l}q^{2k^2+2l^2-2k-2l}q^{4kl}z^{2k}w^{2l}}\over{(q^4;q^4)_l\,(q^4;q^4)_k}}=\cr
&\sum_{l=0}^\infty{{(-1)^{l}q^{2l^2-2l}E_{q^4}(-(q^{2l}z)^2)w^{2l}}\over{(q^4;q^4)_l}}=\sum_{k=0}^\infty{{(-1)^{k}q^{2k^2-2k}E_{q^4}(-(q^{2k}w)^2)z^{2l}}\over{(q^4;q^4)_k}}\cr}$$ which is entire since it is majorized by $E_{q^4}(|w|^2)E_{q^4}(|z|^2)$. Now we wonder whether this expression is lattice-integrable or not. In order to have that we would need that for some $\gamma=(\gamma_1,\,\gamma_2)$ $$\int_{-q\gamma_1\cdot\infty}^{q\gamma_1\cdot\infty}\int_{-\gamma_2\cdot\infty}^{\gamma_2\cdot\infty}(G_{q^2})\cuno(t_1,\,t_2)d_{q^2}\ut$$ is absolutely convergent. For this we would need that 
$$\sum_{h_1=-\infty}^{\infty}\sum_{h_2=-\infty}^\infty q^{2|H|}\Bigl|\sum_{l=0}^\infty{{(-1)^{l}q^{2l^2-2l}E_{q^4}(-(q^{2l+2h_1}q\gamma_1)^2)(q^{2h_2}\gamma_2)^{2l}}\over{(q^4;q^4)_l}}\Bigr|<\infty$$  therefore we have to look at the limit for $h_j\to-\infty$ of the summands. Clearly by the discussion in section 9 of [Koo], we see that we would need to have $\gamma_1=q$. With a similar reasoning we see that $\gamma_2$ must be equal to $1$. Now, for general $z$ and $w$ we have:
$$\eqalign{&(G_{q^2}(\ux)\cuno)(z,\,w)=\sum_{l=0}^\infty{{(-1)^{l}q^{2l^2-2l}w^{2l}(q^{4l}z^2;q^4)}\over{(q^4;q^4)_l}}=\cr
&(z^2;q^4)_{\infty}\sum_{l=0}^\infty{{(-1)^{l}q^{2l^2-2l}w^{2l}}\over{(q^4;q^4)_l(z^2;q^4)_{l}}}=(z^2;q^4)_{\infty}\,_{1}\phi_1(0;z^2;q^4,w^2)\cr}$$ which is the $q^4$ version of the $q$-Bessel function described in [KooSw]. For $(z,\,w)=(q^{2-2r},\,q^{2s})$ with $r\ge0$ and  $s$ any integer, we have, by  the extimates (2.6) and the following extimates in [KooSw],  that $$|(G_{q^2}(\ux)\cuno)(q^{2-2r},\,q^{2s})|=q^{2r(r-1)}q^{4rs}(q^{4r+4};q^4)_{\infty}|\,_{1}\phi_1(0;q^{4r+4};q^4,q^{4r+4s})|$$ For $r\to\infty$ and $s=-r$ this behaves like $q^{2r(r-1)-4r^2}\to\infty$. Hence $G_{q^2}(\ux)$ is not lattice integrable.$\ \spadesuit$\smallbreak
\noindent{\bf Remarks:} An analogue of the symmetry $(2.2)$ for $q$-Bessel functions in [KooSw] holds in our case, namely:
$$E_{q^4}(-x_1^2)\,_1\phi_1(0;x_1^2;q^4,x_2^2)=G_{q^2}(\ux)=E_{q^4}(-x_2^2)E_{q^4}(-x_1^2)=\,_1\phi_1(0;x_2^2;q^4,x_1^2)E_{q^4}(-x_2^2)$$ once we agree that in $\,_1\phi_1$ every time we have a product of type $x_2^l/(x_1^2;q^4)_l$, the terms in $x_1$ have to be taken {\sl before} the terms in $x_2$. Hence  in general one has
$$E_{q}(-x_1-x_2)=E_{q}(-x_2)E_{q}(-x_1)=E_q(-x_1)\,_1\phi_1(0;x_1;q,x_2)=\,_1\phi_1(0;x_2;q,x_1)E_q(-x_2)$$ with the above meaning for $_1\phi_1$ in noncommuting variables. 

\noindent Another equality involving a $_1\phi_1$ and exponentials in $q$-commuting variables is obtained by writing $E_q(-x_1)E_q(-x_2)$ as $$(x_2;q)_\infty\sum_{l=0}^\infty {{(-1)^lq^{{1\over2}(l^2-l)}(q^{-l}x_2;q)_lx_1^l}\over{(q;q)_l}}$$ and using (3.12) in [Koo] with $x_1=-y$ and $x_2=-x$. Then one obtains
$$\sum_{l=0}^\infty{{(-1)^lq^{{1\over2}(l^2-l)}x_1^l(x_2;q)_l}\over{(q;q)_l}}=E_q(x_1x_2)E_q(-x_1)$$ where the sum on the left hand side can be considered as a $\,_1\phi_1$ in noncommuting variables once assumed that $x_1$ always precedes $x_2$ in products. These facts were pointed to me by T. Koornwinder. $\ \spadesuit$\smallbreak
\noindent Of course if a generalized function is $q^2$-integrable then it is lattice integrable for every choice of a lattice. One can easiliy check that for every
$A=(a_1,\ldots,\,a_n)$ in ${\bf R}_{>0}^n$ and every $E=(e_1,\ldots,\,e_n)$  the formal power series $\xe E_{q^4}(-(a_1x_1)^2)\cdots E_{q^4}(-(a_nx_n)^2)$  is lattice-integrable in the $q^2$-lattice generated by $(a_1^{-1}q^{n-1+E^1},\,\ldots,a_j^{-1}q^{n-j+E^j},\,\ldots,\,a_n^{-1})$. Unfortunately, lattice integrability carries a lot of technical work with it whenever one wants to prove anything like translation invariance, for instance.
This is a consequence of the fact that, in order to state that the integral of $\de f(\ux)$ is zero, one needs to keep track of the lattice in which this series is integrable, which in general is not the same as the lattice in which 
$f(\ux)$ is integrable, unless $e_j$ is even for every $j$. Indeed, consider $n=2$ and $f(\ux)=E_{q^4}(-x_1^2)E_{q^4}(-x_2^2)$. Then, $f(\ux)$ is integrable for $(\gamma_1,\,\gamma_2)=(q,\,1)$ while $\partial_2(f(\ux))=-E_{q^4}(-q^2x_1^2){x_2\over{(1-q^2)}}E_{q^4}(-q^4x_2^2)$ is integrable for $(\gamma_1,\,\gamma_2)=(1,\,1)$. However, since ``morally'' the integral of a function which is odd in a variable is zero, we might as well {\sl define} the integral of every odd function to be zero by changing the definition of the integral. Namely, we define the new integral $I'$ to be the integral of the even part of the series $f(\ux)$. We formalize this definition now.\smallbreak
 
\noindent Let $f(\ux)$ be any formal power series in the $x_j$'s. We want to decompose it in $2^n$ series depending on the parity with respect to each variable. Let $$\eqalign{\Pi^{\pm}_j\colon\cove&\to\cove\cr
 f(\ux)&\mapsto{1\over 2}(f(\ux)\pm f(x_1,\ldots,\,x_{j-1},\,-x_j,\,x_{j+1},\ldots,\,x_n))\cr}$$ for every $j$ and for any choice of $\pm$. This makes sense formally, and makes sense even concretely for the series in $C$. Clearly those operators commute, they are projections on the space of power series that are even (resp. odd) in the $j^{\rm th}$ variable, so that  $\Pi^{+}_j\Pi^{-}_j=0$ for every $j$. We define then for every choice of $\beta$ in $\{\pm\}^n$ the operators $\Pi_{\beta}\colon\cov\to\cov$ as $(\Pi^{\beta_1}_1)\circ\cdots\,\circ(\Pi^{\beta_n}_n)$. They are all projections on their image $E_{\beta}$, and clearly the decomposition of the space of power series in the $x_i$'s descends to a decomposition of the space $C$ in $2^n$ spaces that we will call $C_{\beta}$. We also write $V^{\beta}:=C_{\beta}\triangleright 1$. We will denote $\Pi_{(+,\ldots,\,+)}$ by $\Pi_0$ for simplicity.

\noindent In particular $\Pi_0f(\ux)=2^{-n}\sum_{\varepsilon\in\{\pm1\}^n}f(\varepsilon_1x_1,\ldots,\,\varepsilon_n x_n)$ is even in every variable, and we define the integral $I'$ to be the composition $I\circ \Pi_0$.  \smallbreak
\noindent{\bf Remarks:} Since we work in characteristic zero, $I'\cdot f$ is also formally equal to 
$$\int_0^{x_n\cdot\infty}\cdots\int_0^{x_1\cdot\infty}\sum_{\varepsilon\in\{\pm1\}^n}f(\varepsilon_1x_1,\ldots,\,\varepsilon_n x_n)$$
\noindent Clearly the class of $I'$ integrable series in bigger than the class of $I$ integrable series, since all odd series are integrable and their integral is zero. Since $I'$ integrability of a series $f(\ux)$ coincides with $I$ integrability of $\Pi_0f(\ux)$, if $f(\ux)$ is a series which is even in all the variables, then $f(\ux)$ is $I$ integrable $\Leftrightarrow$ $f(\ux)$ is $I'$ integrable since $f(\ux)=\Pi_0f(\ux)$.$\;\spadesuit$\smallbreak
\noindent One can also introduce lattice $I'$ integrability. Again, for series in $C_{(+,\ldots,\,+)}$, lattice $I$ integrability and lattice $I'$ integrability trivially coincide, and for a generic $f(\ux)$, lattice $I'$ integrability trivially coincides with lattice $I$ integrability of $\Pi_0f(\ux)$ in the same lattice.\smallbreak
\noindent We can provide generalizations of Lemma 4-1, Proposition 4-2 and Proposition 4-3 by introducing condition (c') for an $f$ such that $\Pi_0f\in C'$:

\noindent (c') For some $\eta>0$, there exists for each $J=(j_1,\ldots,\,j_n)\in{\bf Z}_{\ge0}^n$ and $\beta\in\{\pm\}^n$  such that  $j_k$ is even (resp. odd) if $\beta_k=+$ (resp. $-$), some constant $C_J$ such that $$|(D_{1,q^2}^{j_1}\cdots D_{n,q^2}^{j_n}(\Pi_{\beta}f)\cuno)(z_1,\ldots,\,z_n)|\le K_J\prod_{k=1}^n(1+|z_k|^2)^{-(1+\eta)}$$ if $\uz\in{\bf R}^n$.\smallbreak 
\noindent Then we have
 
\proclaim Lemma 5-1. Let $f\in C'$ satisfy condition {\rm (c')}. Then for every $J\not=\underline0$ there holds: $(I'(\partial_n^{j_n}\cdots\partial_1^{j_1}f))\triangleright 1\equiv0$, so that we can conclude that $I'(\partial_n^{j_n}\cdots\partial_1^{j_1}f)=0$ .\par

\noindent{\bf Proof: } $I'(\partial_n^{j_n}\cdots\partial_1^{j_1}f)=I(\partial_n^{j_n}\cdots\partial_1^{j_1}\Pi_{\beta}f)$ for $\beta$ related to $J$ as in condition (c'). $\square$

\proclaim Proposition 5-2. Let $f\in C'$ satisfy {\rm (c')}. Then $(id\otimes I')\Delta(f)=1\otimes (I'f)$. $\quad\square$\par

\proclaim Proposition 5-3. Let $f(\ux)\in C$ such that {\rm (c)} holds for every $\eta>0$. Then the statement of Proposition 4-2 holds for every polynomial $p(\ux)$ times $f(\ux)$. $\ \square$\par 

\noindent We also have another invariance property, that is analogous of the classical property (for $n=1$): 
$$\int_{-\infty}^{\infty}{1\over2}(f(x)+f(-x))dx=\int_{-\infty}^{\infty}{1\over2}(f(x+y)+f(x-y))dx$$

\proclaim Proposition 5-4. Let $f(\ux)\in C'$ satisfy {\rm (c')}. Then 
$(id\otimes I)(\Pi_0\otimes id)\Delta(f)=1\otimes(I'f)$. If $f(\ux)$ satisfies condition {\rm (c)} for every $\eta>0$ then the statement is true for every series of the form $\xe f(\ux)$.\par
\noindent{\bf Proof: } The proof uses Taylor's formula with summation only on even $j_k$'s.$\ \square$\smallbreak

\noindent Observe that for even $j_k$'s $$\eqalign{&(I(\partial_n^{j_n}\cdots\partial_1^{j_1}f))\cuno=
(I'(\partial_n^{j_n}\cdots\partial_1^{j_1}f))\cuno=\cr
&=q^{-\sum_kJ^k}\int_{-q^{n-1}z_1\cdot\infty}^{q^{n-1}z_1\cdot\infty}\cdots\int_{-z_n\cdot\infty}^{z_n\cdot\infty}D_{1,q^2}^{j_1}\cdots\,D_{n,q^2}^{j_n}f\cuno(t_1,\ldots,\,t_n)d_{q^2}t_n\cdots d_{q^2}t_1\cr}$$
Hence the proposition above is interesting because it can be proved for lattice integrability with simple changes in the hypothesis and in the proof. This reads as follows.

\noindent Let $\gamma=(\gamma_1,\,\ldots,\,\gamma_n)\in{\bf R}^n$. We define the following spaces:

$$C_{\gamma}=\{f(\ux)\in\cove\ |\ f\cuno|_{\uz=\gamma}{\hbox{ is absolutely convergent }}\}$$ and $C_{q^{2K}\gamma}$ as the space of $f(\ux)\in C_{\gamma}$ such that $f\cuno$ can be continued analytically on a domain containing the $q^2$-lattice generated by $\gamma$. Clearly $C_{\gamma}$ is closed with respect to the multiplication, hence it acts on the space $V_{\gamma}$ of power series in commuting variables $z_1,\ldots,\,z_n$ that are absolutely convergent for $z=\gamma$, hence on a polydisc with polyradius $(|\gamma_1|,\,\ldots,\,|\gamma_n|)$.

\noindent Let $f(\ux)$ be a series in $C_{q^{2K}\gamma}$ for a given $\gamma$. Then it makes sense to investigate $I'(f(\ux))\cuno$ at $z_j=q^{n-j}\gamma_j$ and if this expression is absolutely convergent, then we say that $f(\ux)$ is lattice integrable. Actually we would only need $\Pi_0(f)\in C_{q^{2K}\gamma}$ but since we want to compute integrals of products, we keep the restriction on $f(\ux)$.  

\noindent Consider now the lattice version of condition (c):

\noindent(c'') Let $f(\ux)\in C_{q^{2K}\gamma}$ be such that all Jackson partial derivatives $D_{1,q^2}^{j_1}\cdots D_{n,q^2}^{j_n}f\cuno$ with all $j_k$'s even exist on the lattice $L(\gamma)$, and are such that 
$$|(D_{1,q^2}^{j_1}\cdots D_{n,q^2}^{j_n}(f\cuno))(q^{\pm 2k_1}\gamma_1,\ldots,\,q^{\pm2k_n}\gamma_n)|=O(q^{2(1+\eta)K_-})$$ for $k_i\to \infty$, for some $\eta>0$, where $K_-$ is the sum of the $k_j$'s appearing with the minus sign. \smallbreak
\noindent We introduce the equivalence relation $\sim_{\gamma}$ between two expressions $f(\ux)$ and $g(\ux)$ belonging to $C_{q^{2K}\gamma}$ as follows: $$f(\ux)\sim_{\gamma}g(\ux)\quad\Leftrightarrow\quad f\cuno(\uz)=g\cuno(\uz)\quad \forall \uz\in L(\gamma)$$
\proclaim Proposition 5-5. Let $f(\ux)$ satisfy condition {\rm (c'')} for a given $\gamma$, and let $\gamma'$ denote the $n$-tuple $(q^{n-1}\gamma_1,\,\ldots,\,q^{n-j}\gamma_j,\,\ldots,\,\gamma_n)$. Then 

\noindent{\sl(i) For every $J$ such that every $j_k$ is even, $I'(\partial_n^{j_n}\cdots\partial_1^{j_1}f)\sim_{\gamma'}=0$

\noindent (ii) $(id\otimes I)(\Pi_0\otimes id)\Delta(f)\sim_{\gamma'}1\otimes(I'f)$

\noindent Moreover, if $f(\ux)$ satisfies condition {\rm(c'')} for every $\eta$ and for every $J\in{\bf Z}_{\ge0}^n$ then for every monomial $\xe$

\noindent (iii)  $I'(\partial_n^{2j_n}\cdots\partial_1^{2j_1}(\xe f(\ux)))\sim_{\gamma''}0$ for every $J$, where $\gamma''_j=q^{n-j+E^j}\gamma_j$

\noindent (iv) $(id\otimes I)(\Pi_0\otimes id)\Delta(\xe f(\ux))\sim_{\gamma''}1\otimes I'(\xe f(\ux))$ where $\gamma''$ is as above.}\par

\noindent{\bf Proof :} Statements {\sl (i)} and {\sl (ii)} are clear by the remark after the proof of Proposition 4-4. In order to prove {\sl (iii)} we recall that $$\eqalign{&I'(\partial_n^{2j_n}\cdots\partial_1^{2j_1}(\xe f(\ux)))\cuno=\cr
&q^{-4\sum_kJ^k}\int_{-q^{n-1}z_1\cdot\infty}^{q^{n-1}z_1\cdot\infty}\cdots\int_{-z_n\cdot\infty}^{z_n\cdot\infty}D_{1,q^2}^{2j_1}\cdots D_{n,q^2}^{2j_n}(t_1^{e_1}\cdots t_n^{e_n}f\cuno(q^{E^1}t_1,\,\ldots,\,q^{E^n}t_n))d_{q^2}\ut\cr}$$  hence for $\uz=\gamma''$ this expression converges, and it converges to zero. By invariance under $q^2$ shifts of the Jackson integral we get the statement. Statement (iv) follows from statement (iii). $\square$
\smallbreak
\noindent{\bf Remark: } Observe that in the proof of {\sl (iii)} in Proposition 4-5 the lattice in which we compute the equality depends only on the parity of the $e_j$'s and that it would be enough to be able to keep under control the partial Jackson derivatives of $(P_{\beta}f)\cuno$ with $\beta_j=+$ (resp. $-$) if $e_j$ is even (resp. odd). $\;\spadesuit$\smallbreak
\noindent One may check that $E_{q^4}(-x_1^2)\cdots E_{q^4}(-x_n^2)$ satisfies all the conditions of Proposition 4-5. Computations are left to the reader.\vskip0.2cm

\beginsection6. Lattice order integrability

\noindent We are still left with the problem that the $q^2$-Gaussian $G_{q^2}(\ux)$ is not lattice integrable, even with respect to $I'$. We have to weaken again our condition and introduce the concept {\sl lattice order integrability}. To simplify notation, we use analogues of $I'$ instead of $I$. We first give a simple example of what the procedure is, and the formalize the definition.

\noindent{\bf Example 1} Let $f_1(x_1)$ and $f_2(x_2)$ be even lattice integrable power series for $\gamma_1$ and $\gamma_2$ respectively, if we view them as power series in the one dimensional space. Let $f(\ux)=f_2(x_2)f_1(x_1)$. As we have seen in the Example in the previous Section, this is not necessarily lattice integrable. However, we can consider formally the following construction. We write $f_1(x_1)=\sum_kc_kx_1^k$ and $f_2(x_2)=\sum_k b_kx_2^k$. Then, for $\gamma'=(q\gamma_1,q\gamma_2)$ we have $$\eqalign{&\biggl(\int_{-x_2\cdot\infty}^{x_2\cdot\infty}\int_{-x_1\cdot\infty}^{x_1\cdot\infty}f\biggr)\cuno|_{\uz=\gamma'}=\biggl(\int_{-x_2\cdot\infty}^{x_2\cdot\infty}\int_{-x_1\cdot\infty}^{x_1\cdot\infty}\sum_{k,l}c_kb_lq^{kl}x_1^kx_2^l\biggr)\cuno|_{\uz=\gamma'}=\cr
&\biggl(\int_{-x_2\cdot\infty}^{x_2\cdot\infty}\int_{-x_1\cdot\infty}^{x_1\cdot\infty}\sum_{l}f_1(q^lx_1)b_lx_2^l\biggr)\cuno|_{\uz=\gamma'}\cr}$$ If we formally interchange integration and sum, knowing that from lattice integrability of $f_1$ the result is convergent, we can write
$$\eqalign{&\biggl(\int_{-x_2\cdot\infty}^{x_2\cdot\infty}\int_{-x_1\cdot\infty}^{x_1\cdot\infty}f\biggr)\cuno|_{\uz=\gamma'}=\biggl(\sum_l\int_{-x_2\cdot\infty}^{x_2\cdot\infty}\int_{-q^lx_1\cdot\infty}^{q^lx_1\cdot\infty}f_1(x_1)b_l(q^{-1}x_2)^l\biggr)\cuno|_{\uz=\gamma'}\cr}$$ and since the power series $f_2$ is even, the $q^l$-shift in the integration bound can be neglected. If we  write $I_1(x_1):=\int_{-x_1\cdot\infty}^{x_1\cdot\infty}f_1(x_1)$ we have that the above integral is equal to
$$\eqalign{&\biggl(\sum_l\int_{-x_2\cdot\infty}^{x_2\cdot\infty}I_1(x_1)b_l(q^{-1}x_2)^l\biggr)\cuno|_{\uz=\gamma'}=\cr
&\biggl(I_1(q^{-1}x_1)\sum_l\int_{-x_2\cdot\infty}^{x_2\cdot\infty}b_l(q^{-1}x_2)^l\biggr)\cuno|_{\uz=\gamma'}=\biggl(I_1(q^{-1}x_1)\int_{-x_2\cdot\infty}^{x_2\cdot\infty}f_2(q^{-1}x_2)\biggr)\cuno|_{\uz=\gamma'}
\cr}$$ where we inverted again formally the sum over $l$ and the $q^2$-integral. After applying the action on $1$, and evaluating at $\gamma'$ we get that the result is $q$ times the product of the $q^2$-Jackson integrals of the two powers series evaluated at $\gamma_1$ and $\gamma_2$ respectively. $\spadesuit$\smallbreak 
\noindent We are ready now for the definition of lattice order integrability. What we will do is repeatedly applying a one dimensional integral with respect to a noncommutative variable, say $x_j$. If this expression ``has a meaning'' (i.e. this expression applied to $1$ converges after evaluation at $z_j=\gamma_j$)  then we will identify it with a power series in noncommuting variables, in one variable less, and repeat the procedure. Namely:  
\proclaim Definition 6-1. A formal power series $f(\ux)\in C$ is said to be {\sl lattice order integrable} (l.o. integrable) if there is an ordening  of $1,\,\ldots,\,n$, denoted by the corresponding permutation $\sigma\in S_n$, and an $n$-tuple $\gamma\in{\bf R}_{>0}^n$ such that for every $j\in\{1,\,\ldots,\,n\}$ the expression  $\int_{\sigma(j)}(I_{\sigma(j-1)}\cdots I_{\sigma(1)}f)$ is entire, where $\int_{\sigma(k)}g$ and $I_{\sigma(k)}g$ are defined inductively as follows. For a formal power series $f$ in $\{x_1,\,\ldots,\,x_n\}-\{x_{\sigma(1)},\,\ldots,\,x_{\sigma(k-1)}\}$, 
$\int_{\sigma(k)}f$ is the formal expression in  $\{z_1,\,\ldots,\,z_n\}-\{z_{\sigma(1)},\,\ldots,\,z_{\sigma(k)}\}$
defined as
$$\left(\int_{\sigma(k)}f\right)(\uz):=\left(\int_{-x_{\sigma(k)}\cdot\infty}^{x_{\sigma(k)}\cdot\infty}\Pi_0f\right)\cuno\Bigg|_{z_{\sigma(k)}=\gamma_{\sigma(k)}}$$ 
and if $\int_{\sigma(k)}f$ is entire, $I_{\sigma(k)}f$ will be the unique power series  in the noncommuting indeterminates $\{x_1,\,\ldots,\,x_n\}-\{x_{\sigma(1)},\,\ldots,\,x_{\sigma(k)}\}$ such that $(\int_{\sigma(1)}f)=(I_{\sigma(1)}f)\cuno$. If $f(\ux)$ is l.o. integrable we define the {\sl constant} $I''_{(\sigma,\,\gamma)}f:=\int_{\sigma(n)}(I_{\sigma(n-1)}\cdots I_{\sigma(1)}f)$ to be the {\sl lattice-order integral of $f(\ux)$} associated to the order $\sigma$ and the lattice $L(\gamma)$. \par

\noindent Clearly there is quite a difference between $I$ and $I_{(\sigma,\,\gamma)}''$ since $I$ maps formal power series to formal expressions in $x_1,\,\ldots,\,x_n$ while $I''_{(\sigma,\,\gamma)}$ maps l.o. integrable power series to constants. We will see later what is the relation between the two maps, on the space where they are both defined. We will also see in the examples that even if a power series in l.o. integrable for every order, it could still be not lattice integrable. In this case, we will show the relation between $I''_{(\sigma,\,\gamma)}f$ and $I''_{(\sigma',\,\gamma')}f$. 

\noindent Observe that by definition of $\int_l$, power series that are odd in some variables are automatically defined to be l.o. integrable and that the integral will be zero for every choice of $\sigma$. For this reason, we will only investigate lattice order integrability for even power series. We can state a few results about lattice order integrability.

\proclaim Proposition 6-2. Let $f(\ux)$ be an even element of $\cove$ such that, for some $\tau\in S_n$, and for some power series in one indeterminate $f_1,\,\ldots,\,f_n$ we can write $f(\ux)=f_{\rho(1)}(x_{\rho(1)})\cdots f_{\rho(n)}(x_{\rho(n)})$ where $\rho=\tau^{-1}$. Then $f(\ux)$ is l.o. integrable if and only if each $f_j$, viewed as a power series in commuting variables) is entire and lattice integrable. In this case $f(\ux)$ is lattice order integrable for every order $\sigma$ and a suitable lattice depending on $\sigma$. Moreover, one has $$I''_{(\sigma,\,\gamma)}(f(\ux))=q^{{\rm l}(\sigma)+{\rm l}(\tau)}\prod_{j=1}^n\int_{-\gamma_j\cdot\infty}^{\gamma_j\cdot\infty}(f_j)\cuno(t_j)d_{q^2}t_j$$ where ${\rm l}$ denotes the usual length of a permutation.\par
\noindent{Proof :} $(\Rightarrow)$ Suppose that $f(\ux)$ is as in the hypothesis, and that each $f_j$ is entire, and lattice-integrable for a given $\tilde\gamma_j$. We write $f_j(x_j)=\sum_{k}c_{jk}x_j^k$ for every $j$. For an $n-$tuple $K$ and for $p\in\{1,\,\ldots,\,n\}$, we will also write $$K_{\tau,p}:=\sum_{ j>p\atop \tau(j)<\tau(p)}k_j\quad{\hbox{ and}}\quad K^{\tau,p}:=\sum_{j<p\atop \tau(j)>\tau(p)}k_j$$ We fix a $\sigma$. Then for $\sigma(1)=l$ and $\gamma_l=\tilde\gamma_l$ one has 
$$\eqalign{&\Biggl(\int_lf\Biggr)(z_1,\,\ldots,\,z_{l-1},z_{l+1},\,\ldots,\,z_n)\cr
&=2(1-q^2)\sum_{h=-\infty}^{\infty}q^{2h}\gamma_l\sum_{K}c_{1k_1}\cdots c_{nk_n}q^{-K_l}z_1^{k_1}\cdots z_{l-1}^{k_{l-1}}z_{l+1}^{k_{l+1}}\cdots z_n^{k_n}q^{\sum_jk_jK_{\tau,j}}(q^{2h}\gamma_l)^{k_l}\cr
&=2(1-q^2)\sum_{h=-\infty}^{\infty}q^{2h}\gamma_l\sum_{K'}c_{1k_1}\cdots {\hat c}_{lk_n}\cdots c_{nk_n}(q^{-1}z_1)^{k_1}\cdots(q^{-1}z_{l-1})^{k_{l-1}}z_{l+1}^{k_{l+1}}\cdots z_n^{k_n}\cr
&\times q^{\sum_{j\not=l}k_jK'_{\tau,j}}f_l(q^{2h+K'_{\tau,p}+K'\,^{\tau,p}}\gamma_l)\cr}$$
where $K'$ is the $(n-1)$-tuple obtained by $K$ by deleting $\sigma(1)=l$, and $$K'_{\tau,p}:=\sum_{ l\not=j>p\atop \tau(j)<\tau(p)}k_j\quad{\hbox{ and}}\quad K'\,^{\tau,p}:=\sum_{l\not=j<p\atop \tau(j)>\tau(p)}k_j$$ The last equality holds because $$\eqalign{&\sum_pk_pK_{\tau,p}=k_lK_{\tau,l}+\sum_{p>l}k_pK_{\tau,p}+\sum_{p<l}k_pK_{\tau,p}\cr
&=k_lK_{\tau,l}+\sum_{p>l}k_pK'_{\tau,p}+\sum_{p<l}k_pK'_{\tau,p}+\sum_{p<l\atop\tau(l)<\tau(p)}k_pk_l=k_l(K_{\tau,l}+K^{\tau,l})+\sum_{p\not=l}k_pK'_{\tau,p}\cr}$$ Then one can use the convergence of the $q^2$-Jackson integral of $f_l\,\cuno$, together with the fact that the other $f_k$'s are entire and the fact that $K'_{\tau,p}+K'\,^{\tau,p}$ is an even number because the $f_j$'s are even to show that one can invert the order of summation in the above sum, using dominated convergence. One gets:
$$\eqalign{&\Biggl(\int_lf\Biggr)(z_1,\,\ldots,\,z_{l-1},z_{l+1},\,\ldots,\,z_n)\cr
&=\sum_{K'}q^{\sum_{j\not=l}k_jK'_{\tau,j}}c_{1k_1}\cdots  c_{nk_n}(q^{-1}z_1)^{k_1}\cdots(q^{-1}z_{l-1})^{k_{l-1}}z_{l+1}^{k_{l+1}}\cdots z_n^{k_n}\cr
&\times 2(1-q^2)\sum_{h=-\infty}^{\infty}q^{2h}\gamma_lf_l(q^{2h+K'_{\tau,p}+K'\,^{\tau,p}}\gamma_l)\cr}$$ but this is nothing but
$$\eqalign{&=\sum_{K'}q^{\sum_{j\not=l}k_jK'_{\tau,j}}c_{1k_1}\cdots c_{nk_n}(q^{-1}z_1)^{k_1}\cdots(q^{-1}z_{l-1})^{k_{l-1}}z_{l+1}^{k_{l+1}}\cdots z_n^{k_n}q^{-K'_{\tau,p}-K'\,^{\tau,p}}\cr
&\times\int_{-q^{K'_{\tau,p}+K'\,^{\tau,p}}\gamma_l\cdot\infty}^{q^{K'_{\tau,p}+K'\,^{\tau,p}}\gamma_l\cdot\infty}(f_l)\cuno(t_l)d_{q^2}t_l\cr}$$
The above power series is entire since all the $f_j$'s are, and one finds that
$$I_{\sigma(1)}f=\biggl(\int_{-\gamma_l\cdot\infty}^{\gamma_l\cdot\infty}(f_l)\cuno(t_l)d_{q^2}t_l\biggr)\prod_{j\atop \rho(j)\not=l}f_{\rho(j)}(q^{\eta_{\rho(j)}+\theta_{\rho(j)}}x_{\rho(j)})$$ with $${\eta_k=\cases{-1& if $k<\sigma(1)$\cr
0& if $k>\sigma(1)$\cr}}\qquad{\theta_k=\cases{-1 & if $k<l$ and $\tau(k)>\tau(l)$\cr
1 & if $k>l$ and $\tau(k)<\tau(l)$\cr
0& otherwise\cr}}$$ Therefore we are again in the hypothesis of the Proposition, but then in case $(n-1)$. Since the statement in one dimension is obvious, lattice order integrability is proved, considering the following shifted lattice: for every new step we make, the argument of the $f_j$ that still has to be integrated will be shifted by powers of $q$. If one goes through computations one finds that the exponential of $q$ in the shift of the argument of $f_r$ with $r=\sigma(s)$ is $$-\Sigma_s(\sigma,\tau)=-\left[\#\{j<s\,|\ \sigma(s)<\sigma(j)\}+\#\{j<s\,|\ (\sigma(j)-\sigma(s))(\tau\sigma(j)-\tau\sigma(s))<0\}\right]$$ hence the right lattice to integrate is the one defined by $\gamma_{\sigma(s)}=\tilde\gamma_{\sigma(s)}q^{\Sigma_s(\sigma,\tau)}$. In this setting, the integral will be the product of the $q^2
$-Jackson integrals of the $f_j\,\cuno$'s multiplied by a power of $q$ with exponent 
$$\eqalign{&\sum_{s=1}^n\Sigma_s(\sigma,\tau)={\rm l}(\sigma)+\sum_{s=1}^n\#\{j<s\,|\ (\sigma(j)-\sigma(s))(\tau\sigma(j)-\tau\sigma(s))<0\}\cr
&={\rm l}(\sigma)+{\rm l}(\tau)\cr}$$ since the second term in the sum is equal to the cardinality of
$$\{j,\,s\,|\ j<s\}\cap\bigl(\{j,\,s\,|\;\sigma(j)<\sigma(s),\,\tau\sigma(j)>\tau\sigma(s)\}\cup\{j,\,s\,|\;\sigma(s)<\sigma(j),\,  \tau\sigma(s)>\tau\sigma(j)\}\bigr)$$
For the converse of the statement: one sees that if $f(\ux)$ can be written as a product of one dimensional power series, those series have to be entire, and if there is a $\sigma$ such that $f(\ux)$ is lattice order integrable, this means that $f_{\sigma(1)}$ is lattice integrable on $q^{2k_1}\gamma_{\sigma(1)}$,  and so on, for the following $f_j$'s, with shifted argument. By the $\Rightarrow$ part, we see that lattice order integrability has to hold for every $\sigma'$. $\,\square$.

\noindent{\bf Example 2: } By the above Proposition, the formal power series 

\noindent$f(\ux)=x_n^{e_n}E_{q^4}(-a_n^2x_n^2)\cdots x_1^{e_1}E_{q^4}(-a_1^2x_1^2)$ is l.o. integrable for every $\sigma$ and 

\noindent$\gamma_{\sigma(k)}=a_{\sigma(k)}^{-1}q^{(k-1)+\#\{j<k\,|\;\sigma(k)<\sigma(j)\}}$, since in this case $\tau(k)=n-k+1$.
In particular, $G_{q^2}(q^2\ux)$ and all products of type $G_{q^2}(q^2\ux)\xe$ are l.o. integrable for every choice of the order $\sigma$. 
\noindent One has $$I''_{(\sigma,\,\gamma)}f=\cases{b_{q^2}^n\left[\prod_j(q^2;q^4)_{f_j}\right]q^{2n+2|E|}q^{n\choose2}q^{{\rm l}\,(\sigma)}\prod_{j=1}^na_j^{-1-e_j}&if $e_j=2f_j$ for every $j$\cr
0&otherwise\cr}$$ where
$b_{q^2}=(1-q^2)(q^2, -q^2, -1;\,q^2)_{\infty}$ and the result follows by [Koo].
In particular we observe that the result depends on the choice of $\sigma$ only in a straightforward way and that $L(\gamma)$ does not depend on $E$ but only on $\sigma$ and the $a_j$'s. Therefore it makes sense to consider the relation between $I''_{(\sigma,\,\gamma)}f$ and $I''_{(\sigma,\,\gamma)}(E_{q^4}(-\sum_ka_k^2x_k^2))$. One immediately sees that
if all the $e_j$'s are even $$I''_{(\sigma,\,\gamma)}(f)={{\bigl(\prod_j(q^2;q^4)_{f_j}\bigr)q^{2|E|}}\over{\bigl(\prod_ja_j^{e_j}\bigr)}}I''_{(\sigma,\,\gamma)}(E_{q^4}(-\sum_ka_k^2x_k^2))$$
We say in this case (and whenever an equivalence of integrals $I''_{(\sigma,\gamma)}$ holds, with the same $\sigma$ and $\gamma$ on both sides) that the integral $I(x_n^{e_n}E_{q^4}(-a_n^2x_n^2)\cdots x_1^{e_1}E_{q^4}(-a_1^2x_1^2))$ is {\sl ``weakly equivalent''} to $I(E_{q^4}(-\sum_ka_k^2x_k^2))$. In particular $I(x_n^{e_n}E_{q^4}(-q^2x_n^2)\cdots x_1^{e_1}E_{q^4}(-q^2x_1^2))$ is weakly equivalent to $\prod_j(q^2;q^4)_{f_j}I(G_{q^2}(q^2\ux))$. We also want to point out that the above $f(\ux)$ is an example of the fact that one can have l.o. integrability for every order and still not have lattice-integrability.$\;\spadesuit$\smallbreak

\noindent{\bf Example 3:} If we changed the definition of $\int_l$ for l.o. integrability removing $\Pi_0$ we would no longer have the fact that all constants are l.o. integrable, but then l.o. integrability would strongly depend on the choice of $\sigma$, even in the simplest cases. Take for instance $n=2$, and $f(\ux)=f_1(x_1)f_2(x_2)$ where $f_1(x_1)$ is a series that is not lattice integrable, and $f_2(x_2)$ is lattice-integrable, or even $I-$integrable as a series in one variable, with integral equal to zero. For instance, take $f_1(x_1)=x_1$ and $f_2(x_2)=\partial_2^2e_{q^4}(-x_2^2)$. Then for $\sigma={\rm id}$ the power series $f$ would not be l.o. integrable for any choice of $\gamma$, while for $\sigma=(2\ 1)$ we would have $I''_{(\sigma,\,\gamma)}f=0$ for every $\gamma\in({\bf R_+})^2$.$\;\spadesuit$\smallbreak

\noindent{\bf Properties and Remarks:}

\noindent (a) It is easy to check that if $f(\ux)$ is l.o. integrable for the order $\sigma$ and the lattice $L(\gamma)$, then, for every $n$-tuple of nonzero real numbers $(a_1,\,\ldots,\,a_n)$, the power series $f_A(\ux)=f(a_1x_1,\,\ldots,\,a_nx_n)$ is also l.o. integrable for the same order $\sigma$ and for $\gamma$ replaced by $\tilde\gamma$ where $\tilde\gamma_j=a_j^{-1}\gamma_j$ for every $j$. Then one has equivalence of the numbers $\isigmagamma(f)=(a_1\cdots a_n)I''_{(\sigma,\,\tilde\gamma)}(f_A)$.\smallbreak

\noindent (b) It is also obvious that if $f(\ux)$ is l.o. integrable for $\sigma$ and $\gamma$ then the resulting $\isigmagamma(f)$ is invariant under shifts of each $\gamma_j$ by even powers of $q$.

\noindent (c) If $f$ is as in Proposition 6-1, with $\tau={\rm id}$, and if it is l.o. integrable, then it is lattice integrable.

\noindent (d) In the definition of lattice order integrability the requirement on the $\int_{\sigma(k)}f$'s to be entire for every $k$ can be weakened to analyticity. The weaker version  of the definition is left to the reader. $\;\spadesuit$\smallbreak

\noindent By the discussion above, one can conclude that for well behaved series (by this we mean series satisfying condition (c), (c') etcetera) also the integral $I''_{(\sigma,\,\gamma)}$ is invariant under translation.\smallbreak 

\noindent{\bf Remarks: } The whole construction of lattice order integrability may look artificial, and it may seem to be a definition that is useful only in a noncommutative setting. However, this is not the case. One can define a similar concept of integrability also for power series in commutative variables. In this case, the definition can be made much simpler. 
\proclaim Definition 6-3 (T. Koornwinder). Let $f(z_1,\,\ldots,z_n)=\sum_{e_1,\,\ldots,\,e_n}C_Ez_1^{e_1}\cdots,\,z_n^{e_n}$ 
be a formal power series in commuting variables with coefficients in $\bf R$. The order integral of $f$ with respect to the  order ${\rm id}$ and the lattice $L(\gamma)$, for $\gamma$ an $n$-tuple of real positive numbers, $$\int_{-\gamma_n\cdot\infty}^{\gamma_n\cdot\infty}\left(\cdots\left(\int_{-\gamma_1\cdot\infty}^{\gamma_1\cdot\infty}f(t_1,\ldots,\,t_n)d_{q^2}t_1\right)\cdots\right)d_{q^2}t_n$$ is defined recursively as follows, provided the assumptions in each recursive step are satisfied: $I^0_E:=c_E$ and for $j=1,\ldots,\,n$  $f^j_{e_{j+1},\ldots,e_n}(z_j):=\sum_{e_j}I^{j-1}_{e_j,\ldots,\,e_n}z_j^{e_j}$ is analytic on ${\bf R}+iU$ for some open neighbourhood of $0$ in ${\bf R}$, and $I^j_{e_{j+1},\ldots,\,e_n}:=\int_{-\gamma_j\cdot\infty}^{\gamma_j\cdot\infty}f^j_{e_{j+1},\ldots,\,e_n}(t_j)d_{q^2}t_j$ is absolutely convergent. Then $I^n$ is the required integral. Similarly, if $\sigma$ is a permutation of $\{1,\,\ldots,\,n\}$ we define  
$$\int_{-\gamma_{\sigma(n)}\cdot\infty}^{\gamma_{\sigma(n)}\cdot\infty}\left(\cdots\left(\int_{-\gamma_{\sigma(1)}\cdot\infty}^{\gamma_{\sigma(1)}\cdot\infty}f(t_1,\ldots,\,t_n)d_{q^2}t_{\sigma(1)}\right)\cdots\right)d_{q^2}t_{\sigma(n)}$$ \par
\noindent The definition of order integrability with respect to a lattice in the commutative case is due to T. Koornwinder.\smallbreak
\noindent We have an easy result:
\proclaim Proposition 6-4. Let $f(\ux)=\sum_Ec_E\xe\in C$ be even and lattice order integrable with respect to the order $\sigma$ and the lattice $L(\gamma)$. Then $f\cuno(\uz)$ is integrable with respect to the order $\sigma$ and the lattice $L({\tilde\gamma})$ with ${\tilde\gamma}_{\sigma(l)}=q^{\#\{j<l\,|\,\sigma(l)<\sigma(j)\}}\gamma_{\sigma(l)}$. Moreover $q^{{\rm l}(\sigma)}I^n=I''_{(\sigma,\,\gamma)}f$.\par
\noindent{\bf Proof:} By definition $f_{\cuno}(\uz)=\sum_Ec_E z_1^{e_1}\cdots z_n^{e_n}$. Since by definition of lattice order integrability for $f(\ux)$, $$\eqalign{&\bigl(\int_{\sigma(1)}f\bigr)(z_{\sigma(2)},\ldots,\,z_{\sigma(n)})=\cr
&2(1-q^2)\sum_{k=-\infty}^{\infty}q^{2k}\gamma_{\sigma(1)}\sum_Ec_E(q^{-1}z_1)^{e_1}\cdots(q^{-1}z_{\sigma(1)-1})^{e_{\sigma(1)-1}}(q^{2k}\gamma_{\sigma(1)})^{e_{\sigma(1)}}\cdots z_n^{e_n}\cr}$$ can be written as an {\sl entire power series}  in $z_{\sigma(2)},\cdots,\,z_{\sigma(n)}$, for every $e_{\sigma(2)},\ldots,\,e_{\sigma(n)}$ $I^{\sigma(1)}_{e_{\sigma(2)},\ldots,\,e_{\sigma(n)}}$ is a finite constant and $\sum_{e_{\sigma(2)}}I^{\sigma(1)}_{e_{\sigma(2)},\ldots, e_{\sigma(n)}}z_{\sigma(2)}^{e_{\sigma(2)}}$ converges everywhere. By the fact that $\int_{\sigma(2)}f$ can be expressed as an entire power series,  $\sum_{e_{\sigma(2)}}I^{\sigma(1)}_{e_{\sigma(2)},\ldots, e_{\sigma(n)}}z_{\sigma(2)}^{e_{\sigma(2)}}$ is also integrable with respect to $z_{\sigma(2)}$ in the lattice $L(\tilde\gamma)$ where
$${\tilde\gamma}_{\sigma(2)}=\cases{\gamma_{\sigma(2)}& if $\sigma(2)>\sigma(1)$\cr
q\gamma_{\sigma(2)}& if $\sigma(2)<\sigma(1)$\cr}$$ By repeating the same arguments one obtains the statement. $\square$   

\noindent{\bf Remark: } The reader may wonder whether we could have chosen another realization of $\cove$ and of the integrals of elements of $\cove$ other than $f\triangleright1$. Of course one might consider a different representation, or a different choice of the normal form. The advantages of a representation associated to the choice of a normal form is the fact that it is enough to test operators on $1$ to state an equivalence in $\cov$. The advantage of the particular normal form that we have chosen is based on the fact that $C$ and $C_{\gamma}$ are closed under product, hence we have a map from formal expressions in $x_1\,\ldots,\,x_n$ to formal expressions in the $z_1,\,\ldots,\,z_n$ such that on rather big subspaces it comes exactly from an algebra homomorphism. If we had chosen another normal form,  we could no longer extend the representation $\pi$ of $\cov$ on $V$ to a representation of the subset $S$ of $\cove$ such that  $\pi(S)(1)\subset V$.
 
\noindent Take for instance $n=2$, and the representation of $\cov$ on ${\bf R}[[z,\,w]]$ given by:

\noindent $\pi(x_1)(f(z,\,w))=zf(z,\,q^{-1}w)$ and $\pi(x_2)(f(z,\,w)w)=f(z,\,w)w$. This is the representation associated to the choice of the normal form with $x_2$ preceding $x_1$.

\noindent  Then $a=\sum_{k=0}^\infty x_1^k$ and $b=\sum_{l=0}^\infty x_2^l$  belong to $S$ but $ab$ does not belong to $S$ since $\pi(\sum_{k=0}^\infty x_1^k)(\pi(\sum_{l=0}^\infty x_2^l))(1))=\sum_{k,l=0}^\infty q^{-kl}z^kw^l\not\in V$.

\noindent On the other hand it is true that other realizations, even though they cannot be associated to a representation on $C$ may be more useful to give a meaning to  $I\cdot G_{q^2}(\ux)$ $\;\spadesuit$

\beginsection7. The braided Fourier transform

Now we have all ingredients for the introduction of braided Fourier transforms on a subspace of $\cove$. We introduce two transforms, related to each other by a shift in the arguments and the application of the antipode to one of them. As we already said, the first time that a Fourier transform for this kind of algebras appeared was in [KeMa], from which we took inspiration. A quantum Fourier transform has been studied in the one dimensional case in [Koo]  and in [OR] where also a theory of distributions is developed. The transforms that we present are based on [KeMa], [Koo] and [Ch]. One of the goals of this section is to provide an $n$-dimensional analogues to formulae $(8.19)$, $(8.20)$ and $(8.21)$; hence to Theorem 8.1 in [Koo]. The difference with [Koo] is that  in our version, the algebra $\cove\otimes\vece$ has the braided product $(m\otimes m)(\id\otimes\Phi\otimes\id)$ instead of the ordinary one, although in normal form his formulae and ours for $n=1$ coincide. The difference with [KeMa] lies mainly in the fact that our integral is not bosonic (i.e. it does not have trivial braiding with elements of the algebras $\cove$ and $\vece$). We also apply some minor changes like shifting the argument of the exponential and using the antipode in the definition\footnote{$^{(*)}$}{We do this because the dual of a braided Hopf algebra is the opposite of the dual of a genuine Hopf algebra.}. The use of the antipode appears also in [Ch] where the case of finite dimensional braided Hopf algebras is treated. These transforms behave nicely with respect to a convolution product and with respect to the action of $\vec$ on $\cove$. They also respect various classical properties of the Fourier transform. These facts are developed in [KeMa] and in [Ca].\smallbreak
\noindent We say that an element $f(\ux)$ of $\cove$ is {\sl of class} $\cal R$ if $f(\ux)\,\xe$ is $I'$-integrable for every monomial $\xe$. We say that it is {\sl of class} ${\cal R}_{(\sigma,\,\gamma)}$ if for every monomial $\xe$, the power series $f(\ux)\,\xe$ is lattice order integrable for $\sigma$ and $\gamma$. 

\noindent Again, we do not provide a complete classification of $\cal R$  but we give a class for which this makes sense, which is big enough to reach our goal. Indeed, power series satisfying condition (c) of Section 4 belong to $\cal R$, hence product of $e_{q^2}(-x_j^2)$ and polynomials belong to $\cal R$ provided every $e_{q^2}(-x_j^2)$ appears in the product.\smallbreak

\proclaim Definition 7-1. The braided Fourier transforms $F$ and $F_S$ are defined on the class ${\cal R}$ and they have images in $\cove\otimes\vece$. They are given by  
$$\eqalign{&F:=(I'\otimes{\rm id})(m\otimes{\rm id})\Biggl({\rm id}\otimes \exp\biggl(x\,|\,{{i}\over{(1-q^2)}}\bigl(\partial_1,\,\ldots,\,q^{-j+1}\partial_j,\ldots,\,q^{-n+1}\partial_n\bigr)\biggr)\Biggr)\cr 
&F_S:=(I'\otimes S)(m\otimes{\rm id})\Biggl({\rm id}\otimes \exp\biggl(x\,|\,{{iq^2}\over{(1-q^2)}}\bigl(q^{n-1}\partial_1,\,\ldots,\,q^{n-j}\partial_j,\,\ldots,\,\partial_n)\biggr)\Biggr)\cr}$$\par

\noindent For an $f(\ux)\in {\cal R}\cap C_{\beta}$ one has that: $F(f(\ux)):=(I'\otimes id)\Biggl(f(\ux)\,E_{q^2}\bigl(i\sum_{j=1}^n x_j\otimes q^{-(j-1)}\partial_j\bigr)\Biggr)$ $$\eqalign{&=\sum_{e_1,\,\ldots,\,e_n\atop A(e_j)=\beta_j}I\bigl(f(\ux)\,\xe)\otimes {{i^{|E|}q^{-\sum_jE^j}}\over{\prod_{j=1}^n(q^2;q^2)_{e_j}}}\de\cr}$$ and $F_S(f(\ux)):=(I'\otimes S)\Bigl(f(\ux)E_{q^2}\bigl(iq^2\sum_{j=1}^nx_j\otimes q^{(n-j)}\partial_j\bigr)\Bigr)$ 
$$\eqalign{&=\sum_{e_1,\,\ldots,\,e_n\atop
A(e_j)=\beta_j}I\bigl(f(\ux)\,\xe)\otimes {{(-i)^{|E|}q^{|E|^2-\sum_je_jE_j+|E|+\sum_jE_j}}\over{\prod_{j=1}^n(q^2;q^2)_{e_j}}}\partial_1^{e_1}\cdots,\,\partial_n^{e_n}\cr}$$ where $A(e_j)=+$ (resp. $-$) if $e_j$ is even (resp. odd). Here we used that 
$S(\de)=(-1)^{|E|}q^{|E|^2-|E|-\sum_je_jE_j}\partial_1^{e_1}\cdots\partial_n^{e_n}$
 Therefore it is clear that the second components in the tensor product of $F(f(\ux))$ and of $F_S(f(\ux))$ will also have parity $\beta$. 
In order to provide formulae analogous to $(8.21)$ and $(8.19)$ in [Koo], we need to compute $F_S$ for 

\noindent$M(\ux,\,A)=e_{q^4}(-x_1^2)x_1^{a_1}\cdots\,e_{q^4}(-x_n^2)x_n^{a_n}=x_1^{a_1}\cdots x_n^{a_n}e_{q^4}\bigl(-\sum_j(q^{-A_j}x_j)^2\bigr)$ and for 

\noindent$H(\ux,\,A)=e_{q^4}(-x_1^2){\tilde h}_{a_1}(x_1;q^2)\cdots e_{q^4}(-x_n^2){\tilde h}_{a_n}(x_n;q^2)$ where the ${\tilde h}_{a_j}$'s  are the {\sl discrete $q-$Hermite $II$ polynomials} (see [KoSw] and references in there), that are defined by:$${\tilde h}_l(z;\,q):=z^l\,_{2}\!\phi_1(q^{-n},\,q^{-n+1};\,0;\,q^2,\,-q^2z^{-2})= (q;\,q)_l\,\sum_{k=0}^{[{l\over2}]}{{(-1)^kq^{-2kl+k(2k+1)}z^{l-2k}}\over{(q^2;\,q^2)_k(q;\,q)_{l-2k}}}$$ Both $M(\ux,\,A)$ and $H(\ux,\,A)$  satisfy condition (c) of section 4, so that the transform in defined on both series. We first compute the transform $F_S$ on a generic $f(\ux)$. In order to give a meaning to the transform we  apply the realization map $\pi_\gamma$  sending $g(\ux)$ to $g\cuno(\uz)$, followed by evaluation at $\uz=\gamma$,  to the first component of $F_S(f(\ux))$. By the assumption that $f(\ux)\in{\cal R}$ we know that this is well defined so that $(\pi_{\gamma}\otimes {\rm id})F_S(f(\ux))$ is a genuine  power series in the noncommuting $\partial_j$'s. By the computations in the previous section, one obtains, for an $f(\ux)\in C_{\beta}\cap{\cal R}$:
$$\eqalign{&(\pi_{\gamma}\otimes\id)\bigl(F_S(f(\ux))\bigr)=\sum_{e_1,\,\ldots,\,e_n\atop
A(e_j)=\beta_j} {{(-i)^{|E|}q^{\sum_j(e_j^2+e_j)}}\over{ \prod_{k=1}^n(q^2;\,q^2)_{e_k}}}\times\cr
&\times\left(\int_{-\gamma_1q^{n-1+E_1}\cdot\infty}^{\gamma_1q^{n-1+E_1}\cdot\infty}\cdots \int_{-\gamma_nq^{E_n}\cdot\infty}^{\gamma_nq^{E_n}\cdot\infty}f\cuno(\ut)\,\te d_{q^2}\ut\right)\times \partial_1^{e_1}\cdots\,\partial_n^{e_n} \cr}$$ By invariance of the $q^2-$integral we see that the integration bounds do not depend on $E$ but only on the parity of its components, hence they only depend on $\beta$. In particular, for $f(\ux)=f_1(x_1)\cdots\,f_n(x_n)\in{\cal R}\cap C_\beta$  one has 
$$(\pi_{\gamma}\otimes\id)F_S(f(\ux))=\prod_{k=1}^n\biggl[\sum_{A(e_k)=\beta_k}{{(-i)^{e_k}q^{e_k^2+e_k}\partial_k^{e_k}}\over{(q^2;q^2)_{e_k}}} \biggl(\int_{-q^{B(\beta)_k+n-k}\gamma_k\cdot\infty}^{q^{B(\beta)_k+n-k}\gamma_k\cdot\infty} [(f_k)\cuno(t_k)]t_k^{e_k}d_{q^2}t_k\biggr)\biggr]$$ where the product is taken in {\sl increasing} order and  $B(\beta)=(b(\beta)_1,\ldots,\,b(\beta)_n)$ is the $n$-tuple $\{0,\,1\}^n$ such that the $k^{\rm th}$ entry is $0$ (resp. $1$) if $\beta_k$ is even (resp. odd) and $B(\beta)_k=\sum_{j=1}^{k-1}b(\beta)_j$ as usual. 

\noindent Hence we come to an $n-$dimensional version of formula $(8.21)$ in [Koo]. Let $M(\ux,\,A)$ as above. We remind that in this case $\beta_j=+$ (resp. $-$) if $a_j$ is even (resp. odd). Then we use the one-dimensional case in [Koo] to obtain our result. We can do so because when one chooses the normal forms we have chosen for $\cove$ and $\vece$ computations for each factor look exactly the same as in $(8.21)$ and $(9.8)$ in [Koo]. Indeed, expanding in power series the left hand side of $(8.21)$ and using $q^2$ instead of $q$ one has
$$\eqalign{&c_{q^2}(\gamma)q^{-a^2-a}i^ah_{a}(t;q^2)E_{q^4}(-q^4t^2)=\sum_{k\ge0\atop k+a{\rm \;even}}{{i^kq^{k^2+k}}\over{(q^2;q^2)_k}}\biggl(\int_{-\gamma\cdot\infty}^{\gamma\cdot\infty}x^{a+k}e_{q^4}(-x^2)d_{q^2}x\biggr)t^{k}=\cr
&\sum_{k\ge0\atop k+a{\rm \;even}}{{i^kq^{k^2+k}}\over{(q^2;q^2)_k}}c_{q^2}(\gamma)q^{{-(a+k)^2}\over2}(q^2;q^4)_{{a+k}\over2}t^k\cr}$$

\noindent Using the above formula we obtain:
$$\eqalign{&(\pi_{\gamma}\otimes\id)F_S(M(\ux,\,A))=\cr
&\prod_{j=1}^n\left[c_{q^2}\bigl(q^{n-j+B(\beta)_j}\gamma_j\bigr)\sum_{e_j\ge0\atop
e_j+a_j\;{\rm even}}\left({{(-i)^{e_j}q^{e_j^2+e_j}}\over{(q^2;\,q^2)_{e_j}}}q^{{-(e_j+a_j)^2}\over2}(q^2;\,q^4)_{{e_j+a_j}\over2}\right)\partial_j^{e_j}\right]\cr
&=(-i)^{|A|}\left[\prod_{j=1}^nc_{q^2}(q^{n-j+B(\beta)_j}\gamma_j)\right]q^{\sum_ja_j(1-a_j)}\prod_{k=1}^n\left[E_{q^4}(-q^{4}\partial_k^2)h_{a_k}(\partial_k;\,q^2)\right]\cr}\eqno(7.1)$$ where

\noindent (i) The product is taken in {\sl increasing} order.

\noindent (ii) $E_{q^4}(-z^2)=:G_{q^2}(z)$ is the so-called big $q^2$-Gaussian

\noindent (iii) $h_{l}(z;\,q^2)$ is the {\sl discrete $q^2$-Hermite $I$ polynomial} of degree $l$ (see [KoSw] and references in there) and is defined as $$h_{l}(z;\,q^2):=z^l\,_2\!\phi_0(q^{-2l},q^{2-2l};\,;\,q^4,\,q^{4l-2}z^{-2})=(q^2;\,q^2)_l\,\sum_{k=0}^{[{l\over2}]}{{(-1)^kq^{2k(k-1)}z^{l-2k}}\over{(q^4;\,q^4)_k(q^2;\,q^2)_{l-2k}}}$$

\noindent We observe that the only part of $(\pi_\gamma\otimes\id)F_S(M(\ux,\,A))$ involving the  $\gamma_j$'s is  the coefficient $\left[\prod_{j=1}^nc_{q^2}(q^{n-j+B(\beta)_j}\gamma_j)\right]$ that is also equal to $$\biggl(q^{\sum_jB(\beta)_j}I\bigl(e_{q^4}(-\sum_j(q^{B(\beta)_j}\gamma_j)^2)\bigr)\cuno\biggr)|_{\uz=\gamma}=
\Biggl[\int_{-q^{B(\beta)_n}x_n\cdot\infty}^{q^{B(\beta)_n}x_n\cdot\infty}\cdots\int_{-q^{B(\beta)_1}x_1\cdot\infty}^{q^{B(\beta)_1}x_1\cdot\infty}g_{q^2}\Biggr]\cuno|_{\uz=\gamma}$$ hence the coefficient is a shifted integral of the Gaussian $g_{q^2}(\ux)$ where the shift only depends on the parity of the function $M(\ux,\,A)$, i.e. only on the parity of the $a_j$'s. So we conclude that 
$$\eqalign{&F_S(e_{q^4}(-x_1^2)x_1^{a_1}\cdots\,e_{q^4}(-x_n^2)x_n^{a_n})=\cr
&\Biggl[\int_{-q^{B(\beta)_n}x_n\cdot\infty}^{q^{B(\beta)_n}x_n\cdot\infty}\cdots\int_{-q^{B(\beta)_1}x_1\cdot\infty}^{q^{B(\beta)_1}x_1\cdot\infty}g_{q^2}\Biggr]\otimes(-i)^{|A|}q^{\sum_ja_j(1-a_j)}\prod_{k=1}^n\left[E_{q^4}(-q^{4}\partial_k^2)h_{a_k}(\partial_k;\,q^2)\right]\cr}\eqno(7.2)$$
\noindent The above result gives the analogue of the classical reciprocity between $\{$ Gaussian times a monomial $\}$  and  $\{$ rescaled Gaussian times a Hermite  polynomial $\}$ under the Fourier transform in ${\bf R}^n$. From the above result we derive an analogue of formula $(8.19)$ in [Koo], for $H(\ux,\,A)$ defined above. $H(\ux,\,A)$ is also contained in one of the subspaces $C_{\beta}$ since each ${\tilde h}_{a}(x_j;q^2)$ has constant parity. We obtain:
$$\eqalign{&F_S(e_{q^4}(-x_1^2){\tilde h}_{a_1}(x_1;q^2)\cdots e_{q^4}(-x_n^2){\tilde h}_{a_n}(x_n;q^2))=\Biggl[\int_{-q^{B(\beta)_n}x_n\cdot\infty}^{q^{B(\beta)_n}x_n\cdot\infty}\cdots\int_{-q^{B(\beta)_1}x_1\cdot\infty}^{q^{B(\beta)_1}x_1\cdot\infty}g_{q^2}\Biggr]\otimes\cr
&\otimes\prod_{k=1}^n\left[(q^2;\,q^2)_{a_k}\sum_{s_k=0}^{[{1\over2}a_k]}{{(-i)^{a_k}q^{-a_k^2+a_k}}\over{(q^4;\,q^4)_{s_k}(q^2;\,q^2)_{a_k-2s_k}}}h_{a_k-2s_k}(\partial_k;\,q^2)E_{q^4}(-q^4\partial_k^2)\right]\cr
&=\Biggl[\int_{-q^{B(\beta)_n}x_n\cdot\infty}^{q^{B(\beta)_n}x_n\cdot\infty}\cdots\int_{-q^{B(\beta)_1}x_1\cdot\infty}^{q^{B(\beta)_1}x_1\cdot\infty}g_{q^2}\Biggr](-i)^{|A|}q^{\sum_k(a_k-a_k^2)}\otimes\prod_{j=1}^n\bigl(\partial_j^{a_j}E_{q^4}(-q^{4}\partial_j^2)\bigr)\cr}\eqno(7.3)$$ where the last equality follows from $(8.5)$ in [Koo] and the product is taken in increasing order. After applying the realization map one would get:
$$\eqalign{&(\pi_{\gamma}\otimes\id)F_S(e_{q^4}(-x_1^2){\tilde h}_{a_1}(x_1;q^2)\cdots e_{q^4}(-x_n^2){\tilde h}_{a_n}(x_n;q^2))\cr
&=\left[\prod_{j=1}^nc_{q^2}\bigl(q^{n-j+B(\beta)_j}\gamma_j\bigr)\right](-i)^{|A|}q^{\sum_k(a_k-a_k^2)}\otimes\prod_{j=1}^n\bigl(\partial_j^{a_j}E_{q^4}(-q^{4}\partial_j^2)\bigr)\cr}\eqno(7.4)$$
\noindent We have another analogue of $(8.21)$, which is less helpful, though. It can be obtained with the same techniques as formulae $(5.1)$ and $(5.2)$. It reads as follows:
$$\eqalign{&(\pi_\gamma\otimes\id)F_S\bigl(e_{q^4}\bigl(-\sum_jx_j^2\bigr)x_1^{a_1}\cdots x_n^{a_n}\bigr)=q^{-\sum_ja_jA_j}(\pi_\gamma\otimes\id)F_S(M(q^{A_1}x_1,\ldots,q^{A_n}x_n,A))\cr
&=q^{-\sum_j(a_j+1)A_j}\sum_{e_1,\,\ldots,\,e_n\atop
A(e_j)=\beta_j}{{(-i)^{|E|}q^{|E|^2+|E|}}\over{\prod_j(q^2;q^2)_{e_j}}}\times\cr
&\times\left[(I(M(\ux,\,A))\xe)\cuno(q^{A_1}\gamma_1,\,\ldots,\,q^{A_n}\gamma_n)\right](q^{-A_n}\partial_n)^{e_n}\cdots (q^{-A_1}\partial_1)^{e_n}\cr
&=q^{-\sum A_j(a_j+1)}\biggl(\prod_{j=1}^nc_{q^2}(q^{n-j}\gamma_j)\biggr)(-i)^{|A|}q^{\sum (a_j^2-a_j)}\times\cr
&\prod_{k=1}^n\left[E_{q^4}(-q^{4-2A_k}\partial_k^2)\,h_{a_k}(q^{-A_k}\partial_k;\,q^2)\right]\cr}\eqno(7.5)$$
where again the product is taken in {\sl increasing} order. Hence by the same reasoning as before we can conclude that 
$$\eqalign{&F_S\Bigl(e_{q^4}\bigl(-\sum_jx_j^2\bigr)\,x_1^{a_1}\cdots\,x_n^{a_n}\Bigr)=q^{-\sum_j(a_j+1)A_j}\left[\int_{-x_n\cdot\infty}^{x_n\cdot\infty}\cdots \int_{-x_1\cdot\infty}^{x_1\cdot\infty}g_{q^2}\right]\otimes\cr
&\otimes (-i)^{|A|}q^{\sum_ja_j(1-a_j)}\prod_{k=1}^n\left[E_{q^4}(-q^{4-2A_k}\partial_k^2)h_{a_k}(q^{-A_k}\partial_k;\,q^2)\right]\cr}\eqno(7.6)$$
\noindent This result is less satisfactory because although the coefficients of the $x_j$'s do not depend on the $a_k$'s, the shift in the $\partial_k$'s depends on the $a_k$'s and not only on the position, hence it cannot go through to polynomials. Hence one cannot use this result 
in order to obtain another analogue of  $(8.19)$.\smallbreak
\noindent We observe that for well behaved functions, and for an $n$-tuple $A=(a_1,\,\ldots,\,a_n)$ of nonzero real numbers, the braided Fourier transform of $f(a_1x_1,\,\ldots,\,a_nx_n)$ can be obtained by the braided Fourier transform of $f(\ux)$. More precisely, $(\pi_\gamma\otimes\id)F_S(f(a_1x_1,\,\ldots,\,a_nx_n))=\bigl(\prod_ja_j\bigr)^{-1}(\pi_{\tilde\gamma}\otimes\id)\bigl( F_S(f(\ux))\bigr)(a_1^{-1}\partial_1,\,\ldots,\,a_n^{-1}\partial_n)$ where $\tilde\gamma$ denotes the $n-$tuple obtained by $\gamma$   multiplying each component $\gamma_j$ by $a_j$. Clearly, similar results holds for $(\pi_{\gamma}\otimes\id)F$, hence they hold for $F_S$ and for $F$. This was used to derive formula $(7.5)$.\smallbreak
\noindent Now we want to compute the braided Fourier transform for  monomials or polynomials times a $q^2-$Gaussian of type $G_{q^2}(\ux)$. We cannot use the same definition since $G_{q^2}(\ux)$ is not even lattice integrable. Therefore, we introduce a weaker notion of braided Fourier transform.
\proclaim Definition 7-2. The ``weak'' 
braided Fourier transforms $F''(\sigma,\gamma)$ and $F_S''(\sigma,\gamma)$ are defined on the class ${\cal R}_{(\sigma,\,\gamma)}$. They map this class to $\vece$ and they are defined as
$$\eqalign{&F''(\sigma,\gamma):=(I''_{(\sigma,\,\gamma)}\otimes \id)(m\otimes \id)\Biggl(\id\otimes \exp\biggl(x\,|\,{{i}\over{(1-q^2)}}\bigl(\partial_1,\,\ldots,\,q^{1-j}\partial_j,\,\ldots,\,q^{1-n}\partial_n\bigr)\biggr)\Biggr)\cr
&F''_S(\sigma,\gamma):=(I''_{(\sigma,\,\gamma)}\otimes S)(m\otimes \id)\Biggl(\id\otimes \exp\biggl(x\,|\,{{iq^2}\over{(1-q^2)}}\bigl(q^{n-1}\partial_1,\,\ldots,\,q^{n-j}\partial_j,\,\ldots,\,\partial_n\bigr)\biggr)\Biggr)\cr}$$\par
\noindent We will use this new notion in order to derive an $n$-dimensional version of $(8.20)$ in [Koo]. Namely, we will derive the transform $F''(\sigma,\,\gamma)$ of the formal power series $N(\ux,\,A)=E_{q^4}(-q^4x_n^2)\,x_n^{a_n}\cdots E_{q^4}(-q^4x_1^2)\,x_1^{a_1}$ for given positive integers $a_1,\,\ldots,\,a_n$. One has to compute $I''_{(\sigma,\,\gamma)}(N(\ux,\,A)\,\xe)$ for every $E$ with each $e_j\equiv a_j\pmod 2$, for some fixed $\sigma$ and $\gamma$. This makes sense by the computations in Section 4 because $$\eqalign{&N(\ux,\,A)\,\xe=q^{-\sum_jE^j(e_j+a_j)}E_{q^4}(-q^{4-2E^n}x_n^2)x_n^{e_j+a_j}\cdots E_{q^4}(-q^{4-2E^1}x_1^2)x_1^{e_1+a_1}\cr}$$ is l.o. integrable for every choice of $\sigma$, with $\gamma_{\sigma(k)}=q^{(k-1)+A^{\sigma(k)}+\#\{j<k\;|\;\sigma(j)>\sigma(k)\}}$. In particular, since we showed that the resulting integrals differ only by a factor $q^{l(\sigma)}$, we compute it only for $\sigma={\rm id}$ and $\gamma_k=q^{k-1+A^k}$. Then, denoting by $N'(\ux,\,A)$ the power series obtained by $N(\ux,\,A)$ by multiplying the argument of the $E_{q^4}(-x_j^2)$ by $q^{-E^j}$ for every $j$, one has:
$$\eqalign{&F''(\sigma,\,\gamma)(N(\ux,\,A))\cr
&=\sum_{e_1,\,\ldots,\,e_n\atop e_j+a_j=2h_j}{{i^{|E|}q^{-\sum _jE_j-\sum_jE^j(e_j+a_j)}}\over{\prod_j(q^2;q^2)_{e_j}}}\bigl(I''_{({\rm id},\gamma)}(N'(\ux,\,A+E))\bigr)\partial_n^{e_n}\cdots\partial_1^{e_1}\cr}$$
 So that
$$\eqalign{&F''(\sigma,\,\gamma)(N(\ux,\,A))\cr
&=\sum_{e_1,\,\ldots,\,e_n\atop e_j+a_j=2h_j}{{i^{|E|}q^{-\sum _jE^j-\sum_jE^j(e_j+a_j)+2|E+A|}}\over{\prod_j(q^2;q^2)_{e_j}q^{2|E|+|A|}}}q^{\sum_jE^j(a_j+e_j)}\times\cr
&\times\bigl(\prod_{j=1}^n(q^2;q^4)_{h_j}\bigr)I''_{({\rm id},\,\gamma)}\bigl(E_{q^4}(-\sum_jq^{4-2E^j}x_j^2)\bigr)\partial_n^{e_n}\cdots\partial_1^{e_1}\cr
&=(-1)^{|A|}I''_{({\rm id},\,\tilde\gamma)}(E_{q^4}(-\sum_jq^4x_j^2))\prod_{j=n}^1\left[\sum_{e_j+a_j=2h_j}{{(-i)^{e_j}(q^2;q^4)_{h_j}}\over{\prod_j(q^2;q^2)_{e_j}}}\partial_j^{e_j}\right]\cr
&=q^{\sum_j(a_j^2-a_j)}i^{|A|}b_{q^2}^nq^{n\choose2}\left(\prod_{j=n}^1{\tilde h}_{a_j}(\partial_j;q^2)e_{q^4}(-\partial_j^2)\right)\cr}\eqno(7.7)$$ 
Here $\tilde\gamma$ denotes the $n-$tuple such that $\tilde\gamma_k=q^{A^k}\gamma_k$ and the product is taken in decreasing order. For the last equality we used $(9.15)$ and $(8.20)$ in [Koo].\smallbreak 
\noindent Using the definition of the $h_{a_j}$'s, formula $(7.7)$ above and formula $(8.17)$ in [Koo] one also gets the following result:
$$\eqalign{&F''(\sigma,\,\gamma)(E_{q^4}(-q^4x_n^2)\,h_{a_n}(x_n;q^2)\cdots\,E_{q^4}(-q^4x_1^2)\,h_{a_1}(x_1;q^2))\cr
&=i^{|A|}q^{\sum_j(a_j^2-a_j)}I''_{(\sigma,\,\tilde\gamma)}\bigl(E_{q^4}(-q^4\sum_jx_j^2)\bigr)\prod_{j=n}^1\partial_j^{a_j}e_{q^4}(-\partial_j^2)\cr}\eqno(7.8)$$ where the product is taken in {\sl decreasing} order and $\tilde\gamma$ is  given by $\tilde\gamma_k=q^{A^k}\gamma_k$ for every $k$ as before.\smallbreak
\noindent{\bf Remark:} One could avoid the definition of the weak transform by finding other ways to give a meaning to the integral of power series $f(\ux)$ which is not lattice-integrable. One of the ways could be checking whether, for another choice of the normal form on $\cove$, we could provide another realization $\pi'$ for which both $\pi'(f)$ and $\pi'(I'(f(\ux)\,\xe))$ for every $\xe$ are convergent, or at least converge on a lattice. For instance if we choose the reverse normal form for monomials in $\cov$, then $\pi'N(\ux,\,A)$ would be integrable on a lattice, being the product of $E_{q^4}(-b_j^2z_j^2)z_j^{a_j}$ up to a constant. One could define an equivalence relation on $\vece$  by saying that two elements are the same if for every realization they either both diverge, or they both converge to the same element (at least on a $q^2-$lattice), and there is at least one realization for which they coincide. This is a subject that has to be further developed. These ideas are due to T. Koornwinder.$\,\spadesuit$   
 
\beginsection8. Integral on $\vec$ and Inverse Transform  

We want now to build an inverse for the braided Fourier transforms, at least on a subspace of ${\cal R}$ and ${\cal R}_{(\sigma,\,\gamma)}$. In order to do this one needs an integral on $\vece$, and to develop the theory in this algebra. Since there is a symmetry between $\cove$ and $\vece$, one can simply repeat the definitions and computations keeping in mind that whenever we had a left action involving $\cove$, we will need a right action in the case of $\vece$. We will only provide the necessary formulae, while the properties and the proofs of similar statements as those of Section 4, 5 and 6 are left to the reader. We observe that all the results in this Section can be achieved both by direct computation or by using the symmetry $\psi\colon\cove\to\vece$ defined in Section 2.\smallbreak
\noindent Just as for $\cove$, there are partial derivatives on $\vece$. Those are given on a series $g(\underline\partial)$ by the generalized coefficient of $1\otimes \partial_i$ in $\Delta(g(\upar))$, i.e. for a monomial $\de$, and for $j\in\{1,\ldots,\,n\}$
$$(\de)\leftharpoonup D_i=[e_i]_{q^2}\partial_n^{e_n}\cdots\,\partial_{i+1}^{e_{i+1}}\partial_i^{e_i-1}(q\partial_{i-1})^{e_{i-1}}\cdots\,(q\partial_1)^{e_1}$$
i.e. for $g(\upar)\in\vece$:  $$g(\upar)\leftharpoonup D_i={{g(\partial_n,\ldots,\,q^2\partial_i,\ldots,\, q^2\partial_1)-g(\partial_n,\ldots,\,\partial_i,q^2\partial_{i-1},\ldots,   q^2\partial_1)}\over{(q^2-1)}}\partial_i^{-1}$$
\noindent The operators $D_i$ define a right action of $\cov$ on $\vec$ such that $x_j$ acts as $D_j$. There holds a right version of Taylor's formula, namely:
$$\Delta(g(\upar))=g(\Delta(\upar))=g(\upar)\leftharpoonup\left(\sum_{e_1,\ldots,\,e_n\ge0}{{D_1^{e_1}\cdots\,D_n^{e_n}\otimes\de}\over{[e_1]_{q^2}!\cdots\,[e_n]_{q^2}!}}\right)$$
\noindent A partial inverse for $D_i$ is then given by the (indefinite) $q^2$-integral acting from the right:
$$g(\upar)\leftharpoonup\int_0^{\partial_i}:=(1-q^2)\sum_{k=0}^{\infty}g(\partial_n,\ldots,\,q^{2k}\partial_i,\,q^{-2}\partial_{i-1},\ldots,\,q^{-2}\partial_1)q^{2k}\partial_i$$ and again as in the case of $\cov$, one can define $\int_0^{-\partial_i}$, $\int_0^{q^pa}$. The global integral is formally obtained as the limit for all $r_j\to\infty$ of:
$$g(\upar)\leftharpoonup\int_{-q^{-2r_1}\partial_1}^{q^{-2r_1}\partial_1}\cdots\,\int_{-q^{-2r_n}\partial_n}^{q^{-2r_n}\partial_n}$$ where $$g(\upar)\leftharpoonup\int_{-\partial_j}^{\partial_j}:=\left(g(\upar)\leftharpoonup\int_0^{\partial_j}\right)\,-\,\left(g(\upar)\leftharpoonup\int_0^{-\partial_j}\right)$$
\noindent Hence we have formally
$$\eqalign{&g(\upar)\leftharpoonup\int^{\upar}:=g(\upar)\leftharpoonup\int_{-\partial_1\cdot\infty}^{\partial_1\cdot\infty}\cdots\,\int_{-\partial_n\cdot\infty}^{\partial_n\cdot\infty}\cr
&=  (1-q^2)^n\sum_{\underline\varepsilon\in\{\pm1\}^n}\sum_{k_n=-\infty}^{\infty}\cdots\sum_{k_1=-\infty}^{\infty}g(\varepsilon_nq^{2k_n}\partial_n,\,\ldots,\,\varepsilon_1q^{2k_1}\partial_1)q^{2|K|}\partial_1\cdots\partial_n\cr}$$ 

\noindent As in Section 4, one can define an action of $\vec$ on the power series in the $n$ commuting indeterminates $z_1,\,\ldots,\,z_n$, in order to give a meaning to the integral. One can use this action to define integrability, lattice integrability and lattice order integrability as in Section 5 and 6, and we leave this to the reader. The action will be the right regular action after the choice of a normal form. It is denoted by $\triangleleft$ and it is  defined on monomials as:
$$(z_1^{k_1}\cdots z_n^{k_n})\triangleleft\de:=q^{\sum_jK_je_j}z^{e_1+k_1}\cdots z^{e_n+k_n}_n=z_1^{k_1}(q^{K_1}z_1)^{e_1}\cdots z_n^{k_n}(q^{K_n}z_n)^{e_n}$$ \noindent For brevity we will denote $1\triangleleft g(\upar):=\unoc g(\uz)$  
for every expression $g(\upar)$ for which the action on $1$ makes sense.

\noindent As in Section 5 we construct projections $P^{\pm}_j$ defined on $\vec$ and $\vece$ for every $j=1,\ldots,\,n$ and every choice of $+$ or $-$ as follows:
$$\eqalign{P^{\pm}_j\colon\vece&\to\vece\cr
g(\upar)&\mapsto {1\over2}[g(\upar)\pm g(\ldots,\,-\partial_j,\,\ldots)]\cr}$$ 
\noindent Again, the $P^{\pm}_j$'s commute with each other, they are projections on the subspaces of $\vece$ consisting of  even (resp. odd) elements  with respect to the $j^{\rm th}$ variable, and $P^+_jP^-_j=0$ for every $j$. Then, for every $\beta\in\{-,\,+\}^n$ we define $P_\beta:=P^{\beta_1}_1\circ\cdots\circ P^{\beta_n}_n$. The $P_{\beta}$'s are all projections on their image $G_{\beta}$ and clearly the decomposition $\vece=\bigoplus_{\beta}G_{\beta}$ corresponds to the decomposition of ${\bf C}[[\uz]]$ in series that are either odd or even in each variable after applying the action on $1$. We write $P_0$ for $P_{(+,\,\ldots,\,+)}$. 
\noindent We can define again the integral $J'$ defined as applying $\int^{\upar}$ only to the even part of the power series.

\noindent In particular one can check that for a $g(\upar)$ for which this makes sense one has:
$$\eqalign{&\unoc\bigl(P_0g(\upar)\leftharpoonup J'\bigr)=\unoc\biggl( \Bigl(P_{0}\bigl(g(\upar)\bigr)\Bigr)\leftharpoonup \int^{\upar}\biggr)=\int_{-z_1q^{n-1}\cdot\infty}^{z_1q^{n-1}\cdot\infty}\cdots\int_{-z_nq^{n-n}\cdot\infty}^{z_nq^{n-n}\cdot\infty}\,\unoc(P_0(g))(\ut)d_{q^2}\ut\cr}$$
and 
$$\unoc(g\leftharpoonup D_1^{j_1}\cdots D_n^{j_n})=(D_{q^2}^J(\,\unoc g))(q^{J^1}z_1,\ldots,\,q^{J^n}z_n)$$ 
\noindent We also observe that under the ``symmetry'' $\psi\colon\cove\to\vece$ mapping $x_j$ to $q^{-j+{1\over2}(n-1)}\partial_{n-j+1}$ one has
$$\psi(I\cdot f(x_1,\ldots,x_n))=q^{-n}\biggl[\bigl(\psi(f(\ux))(\upar)\bigr)\leftharpoonup\int_{-q^{-n+1}\partial_1\cdot\infty}^{q^{-n+1}\partial_1\cdot\infty}\cdots\int_{-q^{-n+2j-1}\partial_j\cdot\infty}^{q^{-n+2j-1}\partial_j\cdot\infty}\cdots\int_{-q^{n-1}\partial_n\cdot\infty}^{q^{n-1}\partial_n\cdot\infty}\biggr]\eqno(8.1)$$

\noindent One defines the integral $J''_{(\sigma,\,\gamma)}$ for l.o. integrable power series as the right handed version of $I''_{(\sigma,\,\gamma)}$. We will need $J''_{(\sigma,\,\gamma)}(E_{q^4}(-q^4\partial_1^2)\partial_1^{a_1}\cdots,\,E_{q^4}(-q^4\partial_n^2)\partial_n^{a_n})$. One checks as for the case of $\cove$ that the integrand is actually is l.o. integrable for every $\sigma$ for a suitable choice of $\gamma$, and that the results differ only by a factor $q^{{\rm l}(\sigma)}$. In particular, for $\sigma={\rm id}$ one needs $\gamma_k=q^{(k-1)+A^k}$ and one gets 
$$J''_{(\sigma,\,\gamma)}(E_{q^4}(-q^4\partial_1^2)\partial_1^{a_1}\cdots,\,E_{q^4}(-q^4\partial_n^2)\partial_n^{a_n})=\cases{q^{n\choose2}b_{q^2}^n\prod_j(q^2;q^4)_{b_j}& if $a_j=2b_j$  $\forall\,j$\cr
0& otherwise\cr}$$ hence $$J''_{(\sigma,\,\gamma)}(E_{q^4}(-q^4\partial_1^2)\partial_1^{a_1}\cdots,\,E_{q^4}(-q^4\partial_n^2)\partial_n^{a_n})=\bigl(J''_{(\sigma,\,\gamma)}(E_{q^4}(-q^4\sum_j\partial_j^2)\bigr)\prod_j(q^2;q^2)_{b_j}$$ if $a_j=2b_j$ for every $j$.

\noindent Properties like right invariance of the integral, nullity of the integral of the partial derivative of a power series and so on can be proved as in Section 4.\smallbreak

\noindent Now we are ready to introduce an inversion formula for $F_S$ and $F$ and their weak analogues using the symmetry between $\cove$ and $\vece$ and the results at the end of the previous Section. This will provide an analogue of Theorem 8.1 in [Koo]. The first inversion formula for the braided Fourier transform was due to [KeMa], but we will use a different one and the braiding will not appear in our formulation. As we said our direct transforms are also different from Kempf and Majid's one. The reason why we looked for a different formula is that the element $vol$ in [KeMa] is not necessarily convergent. Moreover, their formula does not necessarily work in our case because our integral is not ``bosonic'', i.e. it has a nontrivial braiding with $\cov$, or $\vec$ for $n\ge2$, as the reader can easily check (see also [Ch] for a few remarks about this property of the integral). Another inversion formula, for the case $n=1$ appears in [OR], where the transforms are  shown to be isomorphisms of topological spaces. We could not extend their results, because we have not yet defined a topology on the noncommutative spaces  we are working with, while for $n=1$ the spaces of functions in $x$ or in $\partial$ can be identified with the classical ones. This is material for future research. \smallbreak
\noindent We say that a power series $g(\upar)$ in $\vece$ is {\sl of class} $\cal S$ if every monomial times $g(\upar)$ is $J'-$integrable. This is possible for instance if $g(\upar)$ satisfies conditions similar to condition (c) of Section 4. We say that $g(\upar)$ is {\sl of class} ${\cal S}_{(\sigma,\,\gamma)}$ if there is an order $\sigma$ and a lattice $L(\gamma)$ such that every monomial  times $g(\upar)$ are lattice order integrable for $\sigma$ and $\gamma$.

\proclaim Definition 8-1. We define the linear maps $G,\,G_S\,\colon {\cal S}\to\cove\otimes\vece$ by
$$\eqalign{&G:=\left[({\rm id}\otimes m) \left(\exp\bigl(x\,|{{i}\over{(1-q^2)}}(\partial_1,\,\ldots,\,q^{-j+1}\partial_j,\,\ldots,\,q^{-n+1}\partial_n)\bigr)\otimes {\rm id}\right)\right]\leftharpoonup\left({\rm id}\otimes J'\right)\cr 
&G_S:=\left[(S\otimes m) \left(\exp\bigl(x\,|\,{{iq^2}\over{(1-q^2)}}(q^{n-1}\partial_1,\,\ldots,\,q^{n-j}\partial_j,\,\ldots,\,\partial_n)\bigr)\otimes {\rm id}\right)\right]\leftharpoonup\left(id\otimes J'\right)\cr}$$\par 
\noindent For instance for $G_S$, the transform of a given $g(\upar)\in{\cal S}$ of fixed parity $\beta$ will be
$$\eqalign{&G_S(g(\upar))=\biggl((S\otimes {\rm id})(E_{q^2}\bigl(iq^2\sum_jq^{(n-j)}x_j\otimes\partial_j)\bigr)g(\upar)\biggr)\leftharpoonup({\rm id}\otimes J')\cr
&=\sum_{e_1,\,\cdots,\,e_n\atop A(e_j)=\beta_j}{{(-i)^{|E|}q^{|E|^2+|E|}q^{\sum_j(n-j)e_j}}\over{\prod_j(q^2;q^2)_{e_j}}}\xe\otimes\left(\de g(\upar)\leftharpoonup({\rm id}\otimes J')\right)\cr}$$

\noindent If $\tau\colon\vece\otimes\cove\to\cove\otimes\vece$ denotes the usual flip operator, putting $c_j=q^{-j+{1\over2}(n-1)}$ for every $j=1,\ldots,\,n$ we observe that for $f(\ux)=f_1(x_1)\cdots f_n(x_n)$ with $f_j$ of parity $\beta_j$, there holds:
$$\eqalign{&\tau(\psi\otimes\psi^{-1})F_S(f(\ux))=\cr
&=\sum_{e_1,\ldots,\,e_n\atop A(e_j)=\beta_j}{{(-i)^{|E|}q^{|E|^2+|E|+\sum_jE_j}}\over{\prod_j(q^2;q^2)_{e_j}}}(c_1^{-1}x_1)^{e_n}\cdots(c_j^{-1}x_j)^{e_{n-j+1}}\cdots(c_n^{-1}x_n)^{e_1}\otimes\cr
&q^{-n}\bigl[f_1(c_1q^{E_1}\partial_n)(c_1\partial_n)^{e_1}\cdots f_n(c_nq^{E_n}\partial_1)(c_n\partial_1)^{e_n}\bigr]\leftharpoonup\int_{-q^{-n+1}\partial_1\cdot\infty}^{q^{-n+1}\partial_1\cdot\infty}\cdots\int_{-q^{n-1}\partial_n\cdot\infty}^{q^{n-1}\partial_n\cdot\infty}\cr}$$ and this is equal to
$$\eqalign{&=\sum_{e_1,\ldots,\,e_n\atop A(e_j)=\beta_j}{{(-i)^{|E|}q^{|E|^2+|E|-\sum_je_jE_j}}\over{\prod_j(q^2;q^2)_{e_j}}}(c_1^{-1}x_1)^{e_n}\cdots(c_j^{-1}x_j)^{e_{n-j+1}}\cdots(c_n^{-1}x_n)^{e_1}\otimes\cr
&\bigl[\partial_n^{e_1}\cdots\partial_1^{e_n}f_1(q^{-E^1}\partial_n)\cdots f_n(q^{-E^n}\partial_1)\bigr]\leftharpoonup\int_{-q^{-n+1+E_n}c_n\partial_1\cdot\infty}^{q^{-n+1+E_n}c_n\partial_1\cdot\infty}\cdots\int_{-c_1q^{n-1+E_1}\partial_n\cdot\infty}^{c_1q^{n-1+E_1}\partial_n\cdot\infty}=\cr
&=\sum_{e_1,\ldots,\,e_n\atop A(e_j)=\beta_j}{{(-i)^{|E|}q^{|E|^2+|E|+\sum_jE^j}}\over{\prod_j(q^2;q^2)_{e_j}}}(c_1^{-1}x_1)^{e_n}\cdots(c_j^{-1}x_j)^{e_{n-j+1}}\cdots(c_n^{-1}x_n)^{e_1}\otimes\cr
&\bigl[\partial_n^{e_1}\cdots\partial_1^{e_n}f_1(\partial_n)\cdots f_n(\partial_1)\bigr]\leftharpoonup\int_{-q^{-n+1+E_n-E_{1}}c_n\partial_1\cdot\infty}^{q^{-n+1+E_n-E^1}c_n\partial_1\cdot\infty}\cdots\int_{-c_1q^{n-1+E_1-E^n}\partial_n\cdot\infty}^{c_1q^{n-1+E_1-E^n}\partial_n\cdot\infty}\cr
\cr}$$ Hence we see that the formal expression of $\tau(\psi\otimes\psi^{-1})F_S(f(\ux))$ coincides with the formal expression of 
$$G_S(f_1(\partial_n)\cdots f_n(\partial_1))((c_1^{-1}x_1,\ldots,c_n^{-1}x_n)\otimes(q^{n-1+E_1-E^1}c_1\partial_n,\ldots,q^{E_n-E^n+1-n}c_n\partial_1))$$ In the same way one shows that the formal expression of $\tau(\psi\otimes\psi^{-1})F(f(\ux))$ coincides with $$G(f_1(\partial_n)\cdots f_n(\partial_1))((c_1^{-1}x_1,\ldots,c_n^{-1}x_n)\otimes(q^{n-1+E_1-E^1}c_1\partial_n,\ldots,q^{E_n-E^n+1-n}c_n\partial_1))$$
We can use the symmetry between $F_S$ and $G_S$, together with formula $(7.2)$ in order to compute $G_S(\partial_n^{a_n}e_{q^4}(-\partial_n^2)\cdots\,\partial_1^{a_1}e_{q^4}(-\partial_1^2))$ for given positive integers $a_1,\,\ldots,\,a_n$ of parity respectively $\beta_1,\ldots,\,\beta_n$. 
Indeed the symmetry tells that:
$$\eqalign{&G_S\bigl(e_{q^4}(-\partial_n^2)\partial_n^{a_n}\cdots e_{q^4}(-\partial_1^2)\partial_1^{a_1}\bigr)=\cr
&(L_{\psi}\otimes L'_{\psi,\beta})\bigl[\tau(\psi\otimes\psi^{-1})\tau F_S(e_{q^4}(-x_1^2)x_1^{a_n}\cdots e_{q^4}(-x_n^2)x_n^{a_1})\bigr]\cr}$$ where $L_{\psi}$ is the shift operator  mapping $x_j$ to $c_jx_j$ and $L'_{\psi,\beta}$ is the shift operator mapping $\partial_j$ to $q^{n-2j+1+B(\beta)^j-B(\beta)_j}c_{n-j+1}^{-1}\partial_j$. Then the above expression is equal to:
$$\eqalign{&(-i)^{|A|}q^{\sum_j(a_j-a_j^2)}\prod_{j=n}^1\bigl[E_{q^4}(-q^4x_j^2)h_{a_j}(x_j;q^2)\bigr]\otimes L'_{\psi,\beta}\psi\int_{-q^{B(\beta)^{1}}x_n\cdot\infty}^{q^{B(\beta)^{1}}x_n\cdot\infty}\cdots\int_{-q^{B(\beta)^{n}}x_1\cdot\infty}^{q^{B(\beta)^{n}}x_1\cdot\infty}g_{q^2}\cr
&=(-i)^{|A|}q^{\sum_j(a_j-a_j^2)}\prod_{j=n}^1\bigl[E_{q^4}(-q^4x_j^2)h_{a_j}(x_j;q^2)\bigr]\otimes\cr
&q^{-n}L'_{\psi,\beta}\biggl(g_{q^2}(c_1\partial_n,\ldots,c_n\partial_1)\leftharpoonup\int_{-q^{B(\beta)^1-n+1}\partial_1\cdot\infty}^{q^{B(\beta)^1-n+1}\partial_1\cdot\infty}\cdots\int_{-q^{B(\beta)^n+n-1}\partial_n\cdot\infty}^{q^{B(\beta)^n+n-1}\partial_n\cdot\infty}\biggr)\cr}$$ So we can conclude that
$$\eqalign{&G_S\bigl(e_{q^4}(-\partial_n^2)\partial_n^{a_n}\cdots e_{q^4}(-\partial_1^2)\partial_1^{a_1}\bigr)=\cr
&(-i)^{|A|}q^{\sum_j(a_j-a_j^2)}\prod_{j=n}^1\bigl[E_{q^4}(-q^4x_j^2)h_{a_j}(x_j;q^2)\bigr]\otimes\biggl(g_{q^2}\leftharpoonup\int_{-q^{B(\beta)_1}\partial_1\cdot\infty}^{q^{B(\beta)_1}\partial_1\cdot\infty}\cdots\int_{-q^{B(\beta)_n}\partial_n\cdot\infty}^{q^{B(\beta)_n}\partial_n\cdot\infty}\biggr)\cr}\eqno(8.2)$$ where the second component of the tensor product clearly depends  only on the parity of the $a_j$'s. Of course formula $(8.2)$ can also be easily obtained by direct computation. Using the definition of the $\tilde h_{m}(z;q^2)$, with the same relation as before between the $a_j$'s and the $\beta_j$'s one obtains:
$$\eqalign{&G_S\bigl({\tilde h}_{a_n}(\partial_n;q^2)e_{q^4}(-\partial_n^2)\cdots\,{\tilde h}_{a_1}(\partial_1;q^2)e_{q^4}(-\partial_1^2)\bigr)\cr
&=(-i)^{|A|}q^{\sum_j(a_j-a_j^2)}\prod_{j=n}^1\bigl[E_{q^4}(-q^4x_j^2)x_j^{a_j}\bigr]\otimes\left(g_{q^2}\leftharpoonup\int_{-q^{B(\beta)_1}\partial_1\cdot\infty}^{q^{B(\beta)_1}\partial_1\cdot\infty}\cdots\int_{-q^{B(\beta)_n}\partial_n\cdot\infty}^{q^{B(\beta)_n}\partial_n\cdot\infty}\right)\cr}\eqno(8.3)$$ which is the $\vece$-version of $(7.3)$. By these results we can conlcude that:
\proclaim Proposition 8-2. Let $\beta=(\beta_1,\,\ldots,\,\beta_n)\in\{\pm1\}^n$, $\sigma$ be any order and $\gamma$ be the $n$-tuple with components  given by $\gamma_{\sigma(k)}=q^{B(\beta)^{\sigma(k)}+k-1+\#\{j<k\,|\,\sigma(j)>\sigma(k)\}}$. Then on power series in $\cove$ of type $E_{q^4}(-q^4x_n^2)p_n(x_n)\cdots E_{q^4}(-q^4x_1^2)p_1(x_1)$ where the $p_i$'s are polynomials of fixed parity $\beta_1,\,\ldots,\,\beta_n$, there holds: 
$$G_S\circ F''(\sigma,\,\gamma)={\rm id}\otimes \bigl(I''_{(\sigma,\,\tilde\gamma)}G_{q^2}(q^2\ux)\bigr)\left(g_{q^2}\leftharpoonup\int_{-q^{B(\beta)_1}\partial_1\cdot\infty}^{q^{B(\beta)_1}\partial_1\cdot\infty}\cdots\int_{-q^{B(\beta)_n}\partial_n\cdot\infty}^{q^{B(\beta)_n}\partial_n\cdot\infty}\right)$$ where $\tilde\gamma$ is such that $\tilde\gamma_k=\gamma_kq^{B(\beta)_k}$ for every $k$. Therefore for power series in $\vece$ of the form  $q_n(\partial_n)e_{q^4}(-\partial_n^2)\cdots\,q_1(\partial_1)e_{q^4}(-\partial_1^2)$ where the $q_j$'s are polynomials of fixed parity $\beta_1,\,\ldots,\,\beta_n$, one has
$$(F''(\sigma,\,\gamma)\otimes{\rm id})G_S={\rm id}\otimes \bigl(I''_{(\sigma,\,\tilde\gamma)}G_{q^2}(q^2\ux)\bigr)\left(g_{q^2}\leftharpoonup\int_{-q^{B(\beta)_1}\partial_1\cdot\infty}^{q^{B(\beta)_1}\partial_1\cdot\infty}\cdots\int_{-q^{B(\beta)_n}\partial_n\cdot\infty}^{q^{B(\beta)_n}\partial_n\cdot\infty}\right)$$ with $\tilde\gamma$ as before.\par
\noindent{\bf Proof:} It follows by $(7.8)$ and $(8.2)$.$\quad\square$\smallbreak
\noindent A slightly more general version of this proposition holds by considering a proper $\gamma$ and $q^2$-Gaussians where the argument is multiplied by a nonzero constant.\smallbreak
\noindent We also want to consider another inverse transform, the weak transform inverting $(7.2)$ and $(7.3)$. 

\proclaim Definition 8-3. The ``weak'' transforms $G''(\sigma,\,\gamma)$ and $G''_S(\sigma,\,\gamma)$ map ${\cal S}_{(\sigma,\,\gamma)}$ to $\cove$ and are defined as:
$$\eqalign{&G''(\sigma,\,\gamma):=\cr
&\left[({\rm id}\otimes m) \left(\exp\bigl(x\,|{{i}\over{(1-q^2)}}(\partial_1,\ldots,q^{-j+1}\partial_j,\ldots,q^{-n+1}\partial_n)\bigr)\otimes {\rm id}\right)\right]\leftharpoonup\left({\rm id}\otimes J''_{(\sigma,\,\gamma)}\right)\cr 
&G''_S(\sigma,\,\gamma):=\cr
&\left[(S\otimes m) \left(\exp\bigl(x\,|{{iq^2}\over{(1-q^2)}}(q^{n-1}\partial_1,\ldots,q^{n-j}\partial_j,\ldots,\partial_n)\bigr)\otimes {\rm id}\right)\right]\leftharpoonup\left({\rm id}\otimes J''_{(\sigma,\,\gamma)}\right)\cr}$$ 

\noindent As in formulae $(7.7)$ and $(7.8)$ one finds, for $\sigma={\rm id}$ and $\gamma_k=q^{A^k+k-1}$ (and similarly for different $\sigma$'s)
$$\eqalign{&G''_{({\rm id},\,\gamma)}\bigl(E_{q^4}(-q^4\partial_1^2)\partial_1^{a_1}\cdots E_{q^4}(-q^4\partial_n^2)\partial_n^{a_n}\bigr)\cr
&=q^{\sum_j(a_j^2-a_j)}i^{|A|}\bigl(\prod_{j=1}^n{\tilde h}_{a_j}(x_j;q^2)e_{q^4}(-x_j^2)\bigr)\bigl[J''_{({\rm id},\,{\tilde\gamma})}\bigl(G_{q^2}(q^2\upar)\bigr)\bigr]\cr}\eqno(8.4)$$ and 
$$\eqalign{&G''_{({\rm id},\,\gamma)}\bigl(E_{q^4}(-q^4\partial_1^2)h_{a_1}(\partial_1)\cdots E_{q^4}(-q^4\partial_n^2)h_{a_n}(\partial_n)\bigr)\cr
&=q^{\sum_j(a_j^2-a_j)}i^{|A|}\bigl(\prod_{j=1}^nx_j^{a_j}e_{q^4}(-x_j^2)\bigr)\bigl[J''_{({\rm id},\,{\tilde\gamma})}\bigl(G_{q^2}(q^2\upar)\bigr)\bigr]\cr}\eqno(8.5)$$ where in both formulae $\tilde\gamma_k=q^{k-1}$ for every $k$ and the product is taken in increasing order. These formulae can be obtained by using the symmetry or by direct computation. One has the second inversion property:
\proclaim Proposition 8-4. Let $\beta=(\beta_1,\,\ldots,\,\beta_n)\in\{\pm1\}^n$, $\sigma$  be any order and $\gamma$ be the $n$-tuple with components  given by $\gamma_{\sigma(k)}=q^{B(\beta)^{\sigma(k)}+k-1+\#\{j<k\,|\,\sigma(j)>\sigma(k)\}}$. Then  on power series in $\vece$ of type $E_{q^4}(-q^4\partial_1^2)p_1(\partial_1)\cdots E_{q^4}(-q^4\partial_n^2)p_n(\partial_n)$ where the $p_j$'s are polynomials of fixed parity $\beta_1,\,\ldots,\,\beta_n$, there holds: 
$$F_S\circ G''(\sigma,\,\gamma)=\biggl[\int_{-q^{B(\beta)_n}x_n\cdot\infty}^{q^{B(\beta)_n}x_n\cdot\infty}\cdots\int_{-q^{B(\beta)_1}x_1\cdot\infty}^{q^{B(\beta)_1}x_1\cdot\infty}g_{q^2}\biggr]\biggl(J''_{(\sigma,\,{\tilde\gamma})}\bigl(G_{q^2}(q^2\upar)\bigr)\biggr)\otimes {\rm id}$$ where $\tilde\gamma_{\sigma(k)}=\gamma_{\sigma(k)}q^{B(\beta)^{\sigma(k)}}$. Therefore for power series in $\cove$ having the form  $q_1(x_1)e_{q^4}(-x_1^2)\cdots q_n(x_n)e_{q^4}(-x_n^2)$ where the $q_j$'s are polynomials of fixed parity $\beta_1,\ldots,\beta_n$, one has
$$({\rm id}\otimes G''(\sigma,\,\gamma))F_S=\biggl[\int_{-q^{B(\beta)_n}x_n\cdot\infty}^{q^{B(\beta)_n}x_n\cdot\infty}\cdots\int_{-q^{B(\beta)_1}x_1\cdot\infty}^{q^{B(\beta)_1}x_1\cdot\infty}g_{q^2}\biggr]\biggl(J''_{(\sigma,\,{\tilde\gamma})}\bigl(G_{q^2}(q^2\upar)\bigr)\biggr)\otimes {\rm id}$$ with $\tilde\gamma$ as before.\par
\noindent{\bf Proof:} It follows by $(7.2)$, $(7.3)$, $(8.4)$ and $(8.5)$.$\quad\square$\smallbreak
\noindent One observes that in this case the Plancherel measure is always a product of integrals of $q^2$-Gaussians, but the integration bounds depend on the parity of the power series. So these transforms could also be seen as sine and cosine transforms (see for this also [KooSw] where the $q$-sine and $q$-cosine transforms in commuting variables are defined).\smallbreak

\noindent{\bf Remark} The break in symmetry appearing in $q^2$-integration is a phenomenon that has recently been observed by [CHMW]. Sometimes this lack of symmetry can be avoided, for instance if the generalized function $f(\ux)$ that we want to integrate (and transform) is lattice integrable in $L(\gamma)$ and $L(q\gamma)$. In this case we could replace $q^2$-integration by $q$-integration since 
$\int_{-\gamma\cdot\infty}^{\gamma\cdot\infty}f\cuno(t)d_qt=\int_{-\gamma\cdot\infty}^{\gamma\cdot\infty}f\cuno(t)d_{q^2}t+\int_{-q\gamma\cdot\infty}^{q\gamma\cdot\infty}f\cuno(t)d_{q^2}t$. The new defined integral will be again invariant under translation. The $q$-integral of a $q^2$-Gaussians $g_{q^2}(\ux)$ times a monomial will be similar to the $q^2$-integral of the same expression. Using the $q$-integral in the definition of $(\pi_{\gamma}\otimes{\rm id})F_S$ will provide  results similar to formula $(7.1)$ but with $\prod_j c_{q^2}(q^{n-j+B(\beta)_j}\gamma_j)$ replaced by $\prod_j(c_{q^2}(\gamma_j)+c_{q^2}(q\gamma_j))$. The result will be therefore independent of the parity of the $a_j$'s. However, this approach cannot be used for $G_{q^2}(\ux)$, since we have seen that there is only one $q^2$-lattice for which $G_{q^2}(\ux)$ is lattice order integrable. $\;\spadesuit$    

\beginsection Acknowledgments

\noindent The author  would like to thank Professor T. H. Koornwinder for many discussions,  and  for his useful comments and remarks.\smallbreak

\beginsection References

\noindent[AST] M. Artin, W. Schelter, and J. Tate {\sl ``Quantum deformation of $GL_n$''} Comm. Pure Appl. Math. {\bf 44} (1991) 879-895

\noindent[Ca] G. Carnovale {\sl``On the braided Fourier transform and the convolution product in the $A_n$ case''} to appear

\noindent[Ch] C. Chryssomalakos {\sl ``Remarks on Quantum Integration''} 
Comm. Math. Phys. {\bf 184} (1997) 1-25

\noindent[CHMW] B. L. Cerchiai, R. Hinterding, J. Madore and J. Wess {\sl ``A Calculus Based on a $q$-deformed Heisenberg Algebra''}, QA/9809160

\noindent[ChZu] C. Chryssomalakos and B. Zumino {\sl``Translations, Integrals and Fourier Transforms in the quantum plane''} in A. Ali, J. Ellis and S. Randjbar-Daemi, eds. ``Salamfestschrift'', proceedings of the ``Conference on highlights of particles and condensed matter physics'' ICTP, Trieste, Italy (1993)

\noindent[Du] M. Durdevic {\sl ``On braided quantum groups''} q-alg 9412003

\noindent[FRT] N. Yu. Reshetikhin, L.A. Takhtadjian and L.D. Faddeev {\sl``Quantization of Lie groups and Lie algebras''} English Transl. Leningrad Math. J.{\bf 1} (1990) 193-225

\noindent[HH] M. Hashimoto and T. Hayashi {\sl ``Quantum multilinear algebra''} T\^ohoku Math. J., {\bf 44} (1992), 471-521

\noindent[HW] R. Hinterding and J. Wess {\sl ``$q$-deformed Hemite Polynomials in $q$-Quantum Mechanics''}, QA/9803050

\noindent[KeMa] A. Kempf and S. Majid {\sl ``Algebraic $q$-integration and Fourier theory on quantum and braided spaces''} J. Math. Phys. {\bf 35}(12) (1994) 
6802-6837

\noindent[Koo] T. H. Koornwinder {\sl ``Special Functions and $q$-commuting variables''} Special Functions, $q$-series and related topics, M.E.H Ismail. D.R. Masson, M. Rahman (eds.),  Fields Institute Communications 14, AMS, (1997) 131-166

\noindent[KooSw] T. H. Koornwinder and R. F. Swarttouw {\sl``On $q$-analogues of the Fourier and Hankel transforms''} Trans. Am. Math. Soc. 333 (1992) 445-461

\noindent[KoSw] R. Koekoek and R. Swarttouw {\sl``The Askey-scheme of hypergeometric orthogonal polynomials and its $q$-analogue''} Revised Version of Report 94-05 Delft University of Technology, Faculty TWI, (1996)

\noindent[LyuMa] V. V. Lyubashenko and S. Majid {\sl``Braided groups and quantum Fourier transform''} J. Alg. 166 506-528 (1994)

\noindent[Ma] S. Majid {\sl ``Foundations of quantum group theory''} Cambridge University Press, (1995)

\noindent[MMNNU] T. Masuda, K. Mimachi, Y. Nakagami, M. Noumi and K. Ueno {\sl``Representations of the Quantum Group $SU_q(2)$ and the Little $q$-Jacobi Polynomials''} J. Funct. Anaysis 99, (1991), 357-387

\noindent[OR] M. Olshanetsky and V. Rogov {\sl``The $q$-Fourier transform of $q$-distributions''} IHES Preprint q-alg 9712055 

\noindent[Ro] L. Rowen {\sl ``Ring Theory''} Vol. I, II Academic press (1988)

\noindent[Sch] M. P. Sch\"utzenberger {\sl``Une interp\'etation de certaines solutions de l'\'equation fonctionelle: $F(x+y)=F(x)F(y)$''} C. R. Acad. Sci. Paris 236 (1953) 352-353

\vfill\eject
\end